\documentclass[12pt]{article}
\usepackage{amssymb}
\begin{document}
\baselineskip=16.3pt
\parskip=3pt
\def\mbb{\mathbb}\def\mb{\mathbb}\def\smb{\mathbb}
\def\SGN{\mbox{sgn}}
\def\BF{(\mbb{F})}\def\BF{}
\def\qed{\hfill \hfill \ifhmode\unskip\nobreak\fi\ifmmode\ifinner
\else\hskip5pt\fi\fi\hbox{\hskip5pt\vrule width4pt height6pt
depth1.5pt\hskip 1 pt}}
\def\WW{}\def\hs{}
\def\ul{\underline}
\def\rb{\raisebox}
\def\rd{\mbox{Rad}}
\def\dis{\displaystyle}
\def\ssc{\scriptscriptstyle}
\def\sc{\scriptstyle}
\def\kn{\mbox{ker}}
\def\psp{{}}
\def\pse{{}}
\def\HOM{\mbox{Hom}^+_{\mb{Z}}(\G,\mbb{F})}
\def\bHOM{\ol{\mbox{Hom}}_{\mb{Z}}(\G,\mbb{F})}
\def\hom{\mbox{Hom}^*_{\mb{Z}}(\G,\mbb{F})}
\def\ptl{\partial}
\def\wh{\widehat}
\def\d{\delta}
\def\D{\Delta}
\def\si{\sigma}
\def\Si{\Sigma}
\def\th{\theta}
\def\a{\alpha}
\def\b{\beta}
\def\G{\Gamma}
\def\g{\gamma}
\def\l{\lambda}
\def\wt{\widetilde}
\def\vf{\varphi}
\def\ad{\mbox{ad}}
\def\stl{\stackrel}
\def\ol{\overline}
\def\es{\varepsilon}
\def\ves{\varepsilon}
\def\vsi{\varsigma}
\def\la{\langle}
\def\ra{\rangle}
\def\vt{\vartheta}
\def\for{\mbox{for}}
\def\rar{\rightarrow}
\def\Rar{\Rightarrow}
\def\Lra{\Leftrightarrow}
\def\der{\mbox{Der$\ssc\,$}}
\def\cl{\centerline}
\def\ni{\noindent}
\def\nl{\newline}
\def\bs{\backslash}
\def\vs{\vspace*}
\def\JJ{{\cal J}}
\def\AA{{\cal A}}
\def\BB{{\cal B}}
\def\CC{{\cal C}}
\def\DD{{\cal D}}
\def\GG{{\cal G}}
\def\HH{{\cal H}}
\def\MM{{\cal M}}
\def\NN{{\cal N}}
\def\PP{{\cal P}}
\def\v#1{\ul{#1}}
\def\N{\mathbb{N}}
\def\Z{\mathbb{Z}}
\def\sZ{\mathbb{Z}}
\def\Q{\mathbb{Q}}
\def\R{\mathbb{R}}
\def\C{\mathbb{C}}
\def\sC{\mathbb{C}}
\def\sF{\mathbb{F}}
\def\F{\mathbb{F}}
\ni{\large\bf Poisson brackets and structure of nongraded
Hamiltonian Lie algebras related to locally-finite derivations}
\par\cl{(to Appear in {\it Canadian J.~Math.})}\vs{2pt}\par\ni \ni{\bf Yucai Su
}
\par\ni
{\small Department of Mathematics, Shanghai Jiaotong University,
Shanghai 200030, P.~R.~China\nl (e-mail: ycsu@sjtu.edu.cn)}
\vs{4pt}\par\ni{\small {\bf Abstract}. Xu introduced a class of
nongraded Hamiltonian Lie algebras. These Lie algebras have a
Poisson bracket structure. In this paper, the isomorphism classes
of these Lie algebras are determined by employing a ``sandwich''
method and by studying some features of these Lie algebras. It is
obtained that two Hamiltonian Lie algebras are isomorphic if and
only if their corresponding Poisson algebras are isomorphic.
Furthermore, the derivation algebras and the second cohomology
groups are determined.} \vs{4pt}\par\ni {\it Mathematical Subject
Classification (1991): 17B40, 17B65} \vs{5pt}\par\ni {\bf1. \
Introduction} \vs{-1pt}\par\ni A Lie algebra $(\AA,[\cdot,\cdot])$
is called to have a {\it Poisson bracket structure} if there
exists a commutative associative algebra structure $(\AA,\cdot)$
such that the compatibility condition holds:
$$
[u, v\cdot w]=[u,v]\cdot w+v\cdot[u,w]\;\;\for\;\;u,v,w\in \AA.
\eqno(1.1)$$
The algebra $(\AA,\cdot,[\cdot,\cdot])$ with two algebraic structures
is also called a {\it Poisson algebra}.
Poisson bracket structures have many applications in areas of mathematics
and physics; they are fundamental algebraic structures on phase spaces in
classical mechanics; they are also the main objects in symplectic geometry
(cf.~[Z]).
\par
Let $\mbb{F}$ be a field of characteristic zero.
A Lie algebra $\AA$ is called {\it graded} if
$\AA=\oplus_{\a\in\G}\AA_{\a}$ is a $\G$-graded $\mbb{F}$-vector space
for some abelian group $\G$ such that
$$
{\rm dim\ssc\,}\AA_{\a}<\infty,\;\;\;\;[\AA_{\a},\AA_{\b}]
\subset\AA_{\a+\b}\;\;\;\for\;\;\;\a,\b\in \G.
\eqno(1.2)$$
A classical Poisson algebra $\PP(\ell)$ is a polynomial algebra
$\AA=\mbb{F}[t_1,t_2,\cdots,t_{2\ell}]$ in $2\ell$ variables with
the Lie bracket
$$
[f,g]=\sum_{i=1}^{\ell}(\ptl_{t_i}(f)\ptl_{t_{\ell+i}}(g)-
\ptl_{t_{\ell+i}}(f)\ptl_{t_i}(g))
\;\;\;\for\;\;\;f,g\in\AA,
\eqno(1.3)$$
where $\ptl_{t_i}$ stands for partial derivative $\ptl\over\ptl{\ssc\,}t_i$.
Define
$$
\PP(\ell)_n= \{t_1^{n_1}t_2^{n_2}\cdots
t^{n_{2\ell}}_{2\ell}\,|\,n_i\in\mbb{N},\,
\sum_{i=1}^{2\ell}n_i=n+2\}\;\;\;\for\;\;\;-2\leq n\in\mb{Z},
\eqno(1.4)$$ then $\PP(\ell)$ is a $\mb{Z}$-graded algebra
$\PP(\ell)=\oplus_{n\in\smb{Z}}\PP(\ell)_n$. When we consider only
its Lie algebra structure, this Lie algebra is denoted by
$\HH(\ell)$. Then $\HH(\ell)$ (or the simple Lie algebra
$[\HH(\ell),\HH(\ell)]/\mbb{F}$) is a classical {\it Lie algebra
of Cartan type H} (also called a {\it Hamiltonian Lie algebra})
[K1,\,K2]. Generalizations of graded Hamiltonian Lie algebras have
been studied in [O,\,OZ].
\par
Nongraded Lie algebras appear naturally in the theory of vertex
algebras and their multi-variable analogues, they play important
roles in mathematical physics. Xu [X2] \mbox{constructed} a family
of in general nongraded Hamiltonian Lie algebras based on certain
derivation-simple algebras and locally finite derivations (we
refer to [SXZ] for the classification of derivation-simple
algebras). In [SX], Xu and the author of this paper determined the
isomorphism classes of Poisson algebras constructed in [X2] (two
Poisson algebras are called isomorphic if there exists an
isomorphism which preserves both associative algebra structure and
Lie algebra structure). However, the structure theory of the
Hamiltonian Lie algebras in general does not seem to be
well-developed. Since the Poisson algebras have two compatible
algebraic structures while the Hamiltonian Lie algebras only have
a Lie algebraic structure, the problem of determination of the
isomorphism classes of Hamiltonian Lie algebras is thus more
complicated, and one can see that some special treatments are
needed in order to \mbox{determine} their isomorphism classes.
\par
In [OZ], Osborn and Zhao determined the isomorphism classes of the
graded Hamiltonian Lie algebras under certain finiteness condition
on the skew-symmetric $\mb{Z}$-bilinear forms $\phi_0$. They used
the ``derivation method'' to determine the isomorphism classes of
the Hamiltonian Lie algebras, mainly, they first determined the
derivation algebras of the Lie algebras in order to obtain their
isomorphism theorem. In this paper, we shall determine the
isomorphism classes of in general nongraded Hamiltonian Lie
algebras $\HH(\v\ell,\G)$, where $\v\ell$ is a 7-tuple of
nonnegative integers and $\G$ is some free abelian group, which
correspond to the Lie algebras in [SX] with the skew-symmetric
$\mb{Z}$-bilinear form $\phi$ being zero and $\ell_4=0$. The
reason we choose $\phi=\ell_4=0$ is that the Hamiltonian Lie
algebras look more natural and more explicit, and are therefore
easier for application, and also they are general enough to cover
already most interesting cases (see \S2). The Hamiltonian Lie
algebras considered in [OZ] in case $\phi_0=0$ are the cases of
the Hamiltonian Lie algebras
$[\HH(\v\ell,\G),\HH(\v\ell,\G)]/\mbb{F}$ with
$\v\ell=(\ell,0,\cdots,0)$.
\par
Unlike the graded case, where the sets of {\it ad}-locally finite
elements and {\it ad}-locally nilpotent elements can be
determined, in the nongraded case, the determination of the sets
of {\it ad}-locally finite elements and {\it ad}-locally nilpotent
elements seems to be un-achievable. Here, we use a ``sandwich''
method to estimate them (see Lemma 3.1). By studying some
important features of the Hamiltonian Lie algebras (Lemma 3.4), we
are able to obtain the isomorphism theorem without the need to
know the structure of their derivation algebras. We obtain
\vs{4pt}\par\ni {\bf Main Theorem}. {\it Two Hamiltonian Lie
algebras are isomorphic if and only if their corresponding Poisson
algebras are isomorphic.}
\par
In Section 2, we shall rewrite the presentations of the
above-mentioned Hamiltonian Lie algebras up to certain obvious
isomorphisms, which we call {\it normalized forms}. Then we shall
prove the main theorem in Section 3. In Section 4, we shall use a
different method from those in [F,OZ] to determine the derivation
algebras of the Hamiltonian Lie algebras. The reason we determine
the derivation algebras after the determination of the isomorphism
classes is that we want to emphasize that the determination of the
isomorphism classes does not depend on the determination of the
derivation algebras. Then in the final section, we shall determine
the second cohomology groups of the Hamiltonian Lie algebras (the
second cohomology groups of the Hamiltonian Lie algebras
considered in [OZ] was determined by Jia [J]).
\par\ni
{\small {\it Acknowledgements.} The author would like to thank
Dr.~Xiaoping Xu for suggesting the investigation of this problem
and for instructions, Professor Kaiming Zhao for helpful
discussions. Part of this research was carried out during the
author's visit to Academy of Mathematics and Systems Sciences,
Chinese Academy of Sciences, he wishes to thank the Academy for
hospitality and support. This research was Supported by a NSF
grant no.~10171064 of China and two research grants from Ministry
of Education of China.} \vs{4pt}\par\ni {\bf2. \ Normalized Forms}
\vs{-1pt}\par\ni Before we present the normalized forms of the
Hamiltonian Lie algebras, to better understand general Hamiltonian
Lie algebras, we first explain how one can generalize the
classical Hamiltonian Lie algebras $\HH(\ell)$ defined in (1.3).
\par
For convenience, we denote
$$
\ol i=i+\ell\;\;\;\for\;\;\;1\le i\le\ell.
\eqno(2.1)$$
The constructional ingredients of the classical Hamiltonian Lie algebra
$\HH(\ell)$ are the pairs
$(\AA,\DD)$ consisting of the polynomial algebra
$$
\AA=\mbb{F}[t_1,t_{\ol1},\cdots,t_\ell,t_{\ol\ell}],
\eqno(2.2)$$
and
a finite dimensional space $\DD={\rm span}\{\ptl_{t_i},\ptl_{t_{\ol i}}
\,|\,1\le i\le\ell\}$ of commuting locally finite derivations.
The derivations $\ptl_{t_i}={\ptl\over\ptl{\ssc\,}t_i}$
are called {\it down-grading operators} by its
obvious meaning for $1\le i\le2\ell$. Then the type of derivation pairs
$\{(\ptl_{t_i},\ptl_{t_{\ol i}})\,|\,1\le i\le\ell\}$ for $\HH(\ell)$ is
$$
(d,d),
\eqno(2.3)$$
where ``d'' stands for down-grading operators.
\par
If we replace the polynomial algebra by the Laurant polynomial
algebra
$$
\AA=\mbb{F}[x_1^{\pm1},x_{\ol1}^{\pm1},\cdots,
x_\ell^{\pm1},x_{\ol\ell}^{\pm1}],
\eqno(2.4)$$
and  rewrite (1.3) as
$$
[f,g]=\sum_{p=1}^{\ell}(x_px_{\ol p})^{-1}
(\ptl^*_p(f)\ptl^*_{\ol p}(g)-
\ptl^*_{\ol p}(f)\ptl^*_p(g))
\;\;\;\for\;\;\;f,g\in\AA,
\eqno(2.5)$$
where $\ptl^*_p$ stands for $x_p{\ptl\over\ptl{\ssc\,}x_p}$ for $1\le p\le
2\ell$, then we obtain a Hamiltonian Lie algebra, denoted by
$\ol{\HH}(\ell)$. Now the
derivations $\ptl_p^*$ are called {\it grading operators} by its obvious
meaning, and the type of derivation pairs
$\{(\ptl^*_p,\ptl^*_{\ol p})\,|\,1\le p\le\ell\}$
for $\ol{\HH}(\ell)$ is then
$$
(g,g),
\eqno(2.6)$$
where ``g'' stands for grading operators.
\par
Furthermore, we can replace $\AA$ by a semigroup algebra which is
the tensor product of a Laurant polynomial algebra (2.4) and a
polynomial algebra (2.2):
$$
\AA=\mbb{F}[x_1^{\pm1},t_1,x_{\ol1}^{\pm1},t_{\ol1},
\cdots,x_{\ell}^{\pm1},t_{\ell},x_{\ol\ell}^{\pm1},t_{\ol\ell}],
\eqno(2.7)$$
and replace $\ptl^*_p$ by $\ptl_p=\ptl^*_p+\ptl_{t_p}$
for $1\le p\le2\ell$,
then (2.5) defines a Hamiltonian Lie algebra, denoted by $\wh{\HH}(\ell)$.
The derivation $\ptl_p$ are called {\it mixed operator}, and the type
of derivation pairs
$\{(\ptl_p,\ptl_{\ol p})\,|\,1\le p\le\ell\}$
for $\wh{\HH}(\ell)$ is now
$$
(m,m),
\eqno(2.8)$$
where ``m'' stands for mixed operators.
\par
In the examples above, we can generally denote a monomial as
$$
x^{\a,\v i}=x_1^{\a_1}x_{\ol1}^{\a_{\ol1}}
\cdots x_\ell^{\a_\ell}x_{\ol\ell}^{\a_{\ol\ell}}
t_1^{i_1}t_{\ol1}^{i_{\ol1}}\cdots t_\ell^{i_\ell}
t_{\ol\ell}^{i_{\ol\ell}},
\eqno(2.9)$$
for
$$\a=(\a_1,\a_{\ol1},\cdots,\a_\ell,\a_{\ol\ell})\in\G,\;\;\;
\v i=(i_1,i_{\ol1},\cdots,i_\ell,i_{\ol\ell})\in\JJ,
\eqno(2.10)$$
where $\G$ is an additive subgroup of $\mbb{F}^{2\ell}$ such that
$\G=\{0\}$ in the case of $\HH(\ell)$
(where there are no nonzero grading operators), and $\G=\mb{Z}^{2\ell}$ in
the cases of $\ol{\HH}(\ell)$ and $\wh{\HH}(\ell)$
(where there are nonzero grading operators), and where $\JJ$ is some
semi-subgroup of $\mbb{N}^{2\ell}$ such that $\JJ=\mbb{N}^{2\ell}$
in the cases of $\HH(\ell)$ and $\wh{\HH}(\ell)$
(where there are nonzero down-grading operators), and
$\JJ=\{0\}$ in the case of $\ol{\HH}(\ell)$ (where there are no nonzero
down-grading operators). In all three cases, we can
define operators
$\ptl^*_p=x_p{\ptl\over\ptl{\ssc\,}x_p},\,\ptl_{t_p}=
{\ptl\over\ptl{\ssc\,}t_p}$ and $\ptl_p=\ptl^*_p+\ptl_{t_p}$ such that
$\ptl^*_p=0$ in the case of $\HH(\ell)$ and $\ptl_{t_p}=0$ in the case
of $\ol{\HH}(\ell)$.
\par
With the above examples in mind, we can now give generalizations
of the Hamiltonian Lie algebras as follows.
\par
First for convenience, for $m,n\in\mb{Z}$, we denote
$$
\ol{m,n}=\left\{\matrix{\{m,m+1,\cdots,n\}&\mbox{if \ \ \ }m\le n
\vs{4pt}\hfill\cr
\emptyset\hfill&\mbox{otherwise}.
\hfill\cr}\right.
\eqno(2.11)$$
\par
We shall construct a semigroup algebra $\mbb{F}[\G\times\JJ]$
(cf.~(2.7)), where $\G$ is some free abelian subgroup of an
$\mbb{F}$-vector space $\mbb{F}^n$ and $\JJ$ is some semi-subgroup
of $\mbb{N}^n$, and construct 7 groups of derivation pairs
$\{(\ptl_p,\ptl_{\ol p})\,|\,p\in I_i\}$ for $i\in\ol{1,7}$, where
$I_i$ are some indexing sets such that if we denote each type of
derivation pairs $\{(\ptl_p,\ptl_{\ol p})\,|\,p\in I_i\}$ by
$(T_i,T_{\ol i})$ for $i\in\ol{1,7}$, then the types of derivation
pairs in the order of the groups $\{(\ptl_p,\ptl_{\ol p})\,|\,p\in
I_i\}$ for $i\in\ol{1,7}$ are
$$
\matrix{
(T_1,T_{\ol 1})=(g,g),\;\;\;\;
(T_2,T_{\ol 2})=(m,g),\;\;\;\;
(T_3,T_{\ol 3})=(m,g),\;\;\;\;
(T_4,T_{\ol 4})=(m,m),
\vs{4pt}\hfill\cr
(T_5,T_{\ol 5})=(g,d),\;\;\;\;
(T_6,T_{\ol 7})=(m,d),\;\;\;\;
(T_7,T_{\ol 7})=(d,d).
\hfill\cr}
\eqno(2.12)$$
Then we shall see that (2.3), (2.6) and (2.8) correspond respectively
to the three special cases: (i) $I_7=\ol{1,\ell}$ and
$I_i=\emptyset$ if $i\ne7$, (ii) $I_1=\ol{1,\ell}$ and
$I_i=\emptyset$ if $i\ne1$, and (iii) $I_4=\ol{1,\ell}$ and
$I_i=\emptyset$ if $i\ne4$.
\par
To construct, we let
$$
\v\ell=(\ell_1,\cdots,\ell_7)\in \mbb{N}^7\bs\{0\}.
\eqno(2.13)$$
Set
$$
\iota_0=0,\;\;\;\iota_i=\ell_1+\ell_2+...+\ell_i,\;\;i\in\ol{1,7},
\eqno(2.14)$$
$$
I_{i,j}=\ol{\iota_{i-1}+1,\iota_j}\;\;\;\for\;\;\;i,j\in\ol{1,7},\,\;i\le j.
\eqno(2.15)$$
Denote
$$
I_i=I_{i,i},\;\;\;I=I_{1,7},\;\;\;J=\ol{1,2\iota_7}.
\eqno(2.16)$$
Define the map \rb{6pt}{$\ol{\ \sc\,}\sc\,$}: $J\rar J$ by
$$
\ol p=\left\{\matrix{
p+\iota_7\hfill&\mbox{if \ \ \ }p\in\ol{1,\iota_7},
\vs{4pt}\hfill\cr
p-\iota_7\hfill&\mbox{if \ \ \ }p\in\ol{\iota_7+1,2\iota_7},
\hfill\cr}\right.
\eqno(2.17)$$
(cf.~(2.1)).
For any subset $K$ of $\ol{1,2\iota_7}$, we denote
$$
\ol K=\{\ol p\,|\,p\in K\}.
\eqno(2.18)$$
In particular, we have $J=I\cup\ol I.$
Set
$$
J_i=I_i\cup\ol I_i,\;\;\;J_{i,j}=I_{i,j}\cup\ol I_{i,j}\;\;\;\for\;\;\;
i,j\in\ol{1,7},\,\;i\le j.
\eqno(2.19)$$
Let $\mbb{F}$ be a field of characteristic zero.
We write an element $\a$ of $\mbb{F}^{2\iota_7}$ in the form
$$
\a=(\a_1,\a_{_{\sc \ol1}},\cdots,\a_{\iota_7},\a_{_{\sc
\ol{\iota_7}}}) \mbox{ \ \ \ with \ \ \ }\a_p\in\mbb{F},
\eqno(2.20)$$ (cf.~(2.10)). Set
$$
\es_p=(\d_{1,p},\d_{\ol
1,p},\cdots,\d_{\iota_7,p},\d_{\ol{\iota_7},p})\in
\mbb{F}^{2\iota_7}\;\;\;\for\;\;\;p\in J. \eqno(2.21)$$ For
$\a\in\mbb{F}^{2\iota_7}$ and $K\subset J$, we use $\a_{_K}$ to
denote the vector in $\mbb{F}^{|K|}$ (where $|K|$ is the size of
$K$), obtained from $\a$ by deleting all the coordinate $\a_p$
with $p\in J\setminus K$; for instance,
$$
\a_{_{\sc \{1,3\}}}=(\a_1,\a_3)\in\mbb{F}^2,\;\;\; \a_{_{\sc
\{1,\bar{2},\bar{3}\}}}=(\a_1,\a_{\ol2},\a_{\ol3}) \in\mbb{F}^3.
\eqno(2.22)$$ Sometimes, when the context is clear, we also use
$\a_{_K}$ to denote the vector in $\mbb{F}^{2\iota_7}$ by putting
its $p$th coordinate to be zero for $p\in J\bs K$.
\par
We fix a set $\{\si_p\,|\,p\in J\}$ of
elements in $\mbb{F}^{2\iota_7}$ as follows:
$$
\si_p=\left\{\begin{array}{ll} \es_p+\es_{\ol p}&\mbox{if \ \ \
}p\in I_1\cup I_{3,4},
\vs{4pt}\\
\es_p&\mbox{if \ \ \ }p\in I_2,
\vs{4pt}\\
0 &\mbox{if \ \ \ }p\in I_{5,7},\end{array}\right. \eqno(2.23)$$
and $\si_{\ol p}=\si_p$.
Using the notations (2.9) and (2.23),
the factor $(x_px_{\ol p})^{-1}$ appears
in (2.5) is simply $x^{-\si_p}$ if $p\in I_1$.
If we re-denote $x_p^{-1}$ by $x_p$ (and $x_p$ by $x_p^{-1}$), then
the factor $(x_px_{\ol p})^{-1}$ in (2.5) can be written as
$$
(x_px_{\ol p})^{-1}=x^{\si_p}.
\eqno(2.24)$$
\par
Now we take an additive subgroup $\G$ of
$\mbb{F}^{2\iota_7}$ such that
$$
\a_{_{\sc \ol I_{5,6}\cup J_7}}=0\;\;\;\for\;\;\;\a\in\G,
\eqno(2.25)$$ (this condition is necessary since we require that
$T_{\ol5}=T_{\ol6}=T_7=T_{\ol7}=d$ by (2.12), which means that
$\ptl^*_p=0$, i.e., we shall have $\a_p=0$ if $p\in \ol
I_{5,6}\cup J_7$ for $\a\in\G$ (cf.~(2.2) and
(2.3)$\ssc\,$)$\ssc\,$), and we shall also require that
$$
\si_p\in\G,\ \ \es_q\in\G,\ \ \mbb{F}\es_r\cap\G\ne\{0\}
\;\;\;\for\;\;\;p\in I_{1,4},\,q\in I_{5,6},\,r\in J_{1,4},
\eqno(2.26)$$
where the first condition is necessary since we require that
$x^{\si_p}$ will appear as a factor in the Lie bracket (cf.~(2.5) and (2.24),
also see (2.36)), and where the last two
conditions are called the {\it distinguishable conditions} among the
derivations $\ptl_p$ defined later
in (2.33), which are necessary in order to
guarantee the simplicity of the Hamiltonian Lie algebras (cf.~[X2]).
\par
Note that $\mbb{N}^{2\iota_7}$ is an additive semi-subgroup of
$\mbb{F}^{2\iota_7}$. We take
$$
\JJ=\{\v i=(i_1,i_{\ol1},\cdots,i_{\iota_7}, i_{_{\sc
\ol{\iota_7}}}) \in\mbb{N}^{2\iota_7}\,|\,\v i_{J_1\cup\ol
I_{2,3}\cup I_5}=0\}, \eqno(2.27)$$ (cf.~(2.10)), where the
condition $\v i_{J_1\cup\ol I_{2,3}\cup I_5}=0$ is necessary since
$T_1=T_{\ol1}=T_{\ol2}=T_{\ol3}=T_5=g$ by (2.12), which means that
$\ptl_{t_p}=0$, i.e., we shall have $i_p=0$ if $p\in J_1\cup\ol
I_{2,3}\cup I_5$ for $\v i\in\JJ$ (cf.~(2.4) and (2.6)).
\par
Now we let $\AA=\mbb{F}[\G\times\JJ]$ be the semigroup algebra with basis
$$
\{x^{\a,\v i}\,|\,(\a,\v i)\in\G\times\JJ\},
\eqno(2.28)$$
(cf.~(2.9)),
and the multiplication
$$
x^{\a,\v i}\cdot x^{\b,\v j}=x^{\a+\b,\v i+\v j}\;\;\;
\for\;\;\;(\a,\v i),\,(\b,\v j)\in\G\times \JJ.
\eqno(2.29)$$
Then $\AA$ forms a commutative associative algebra with $1=x^{0,0}$
as the identity element.
Set
$$
\AA_{\a}={\rm span}\{x^{\a,\v i}\,|\,\v i\in\JJ\}\;\;\;\for\;\;\;
\a\in\G.
\eqno(2.30)$$
Then $\AA$ is $\G$-graded $\AA=\oplus_{\a\in\G}\AA_\a$
(but in general $\AA_\a$ is infinite dimensional).
For convenience, we denote
$$
x^{\a}=x^{\a,0},\;\;\;t^{\v i}=x^{0,\v i},\;\;\;t_p=t^{\es_p},\;\;\;
\for\;\;\;\a\in\G,\;\v i\in\JJ,\;p\in J.
\eqno(2.31)$$
In particular,
$$
t^{\v i}=\prod_{p\in J}t_p^{i_p},\;\;\;
x^{\a,\v i}=x^\a t^{\v i},\;\;\;\for\;\;\;\a\in\G,\;\v i\in\JJ,
\eqno(2.32)$$
(cf.~(2.9)).
Define the derivations $\{\ptl_p,\ptl^*_p,\ptl_{t_p}\,|\,p\in J\}$ of $\AA$
by
$$
\ptl_p=\ptl^*_p+\ptl_{t_p}\mbox{ \ \ \ and \ \ \ }
\ptl^*_p(x^{\a,\v i})=\a_p x^{\a,\v i},\ \ \
\ptl_{t_p}(x^{\a,\v i})=i_p x^{\a,\v i-\es_p},
\eqno(2.33)$$
for $p\in J,\,(\a,\v i)\in\G\times\JJ,$
where we treat
$$
x^{\a,\v i}=0\mbox{ \ \ if \ \ }(\a,\v i)\notin\G\times\JJ.
\eqno(2.34)$$
In particular,
$$
\ptl^*_p=0,\;\;\; \ptl_{t_q}=0\;\;\;\for\;\;\;p\in
\ol I_{5,6}\cup J_7,\;\;q\in J_1\cup\ol I_{2,3}\cup I_5,
\eqno(2.35)$$
by (2.25) and (2.27) (cf.~(2.12)).
We call the nonzero derivations $\ptl^*_p$
{\it grading operators}, the nonzero derivations $\ptl_{t_q}$
{\it down-grading operators}, and the derivations $\ptl_r^{\ast}+\ptl_{t_r}$
{\it mixed operators} if both $\ptl_r^{\ast}$ and $\ptl_{t_r}$ are not zero.
Then the types of derivation pairs in the order of the groups
$\{(\ptl_p,\ptl_{\ol p})\,|\,p\in I_i\}$ for $i\in\ol{1,7}$ are
shown as in (2.12).
\par
Now we define the following Lie bracket on $\AA$:
$$
[u,v]=\sum_{p\in I}x^{\si_p} (\ptl_p(u)\ptl_{_{\sc \ol
p}}(v)-\ptl_{_{\sc \ol p}}(u)\ptl_p(v)), \eqno(2.36)$$ for
$u\in\AA_{\a},v\in\AA_{\b}$ (cf.~(2.30), (2.5) and (2.24)), where
$x^{\si_p}$ appears just as in (2.5) and (2.24). Then
$(\AA,[\cdot,\cdot])$ forms a Hamiltonian Lie algebra, denoted by
$\HH(\v\ell,\G)$, and $(\AA,\cdot,[\cdot,\cdot])$ forms a Poisson
algebra. Then $\HH(\v\ell,\G)$ is the normalized form of a class
of in general nongraded Hamiltonian Lie algebra constructed in
[X2]. {}From this definition, one sees that the classical
Hamiltonian Lie algebra $\HH(\ell)$ is simply the Lie algebra
$\HH(\v{\ell'},0)$ with $\v{\ell'}=(0,\cdots,0,\ell)$, and the
Hamiltonian Lie algebras $\ol{\HH}(\ell)$ and $\wh{\HH}(\ell)$ are
respectively $\HH(\v{\ell''},\mb{Z}^\ell)$ and
$\HH(\v{\ell'''},\mb{Z}^\ell)$, where $\v{\ell''}=
(\ell,0,\cdots,0)$, and $\v{\ell'''}= (0,0,0,\ell,0,0,0)$
(cf.~(2.12) and the statement after it). The Hamiltonian Lie
algebras considered in [OZ] in case $\phi_0=0$ are the cases of
the Hamiltonian Lie algebras
$[\HH(\v\ell,\G),\HH(\v\ell,\G)]/\mbb{F}$ with
$\v\ell=(\ell,0,\cdots,0)$.
\par
The Hamiltonian Lie algebras $\HH(\v\ell,\G)$ can also be
viewed
as generalizations of the Lie algebras in [DZ,\,X1,\,Zh]
in the sense that they have some common features stated in Lemma 3.4.
\par
The following theorem was proved in [X2].
\par\ni
{\bf Theorem 2.1}. {\it
The Lie algebra $\HH(\v\ell,\G)$ is central simple, i.e.,
$[\HH(\v\ell,\G),\HH(\v\ell,\G)]/\mbb{F}$
(the derived algebra modulo its center) is simple.}
\vs{5pt}\par\ni
{\bf3. \ Isomorphism Classes}
\vs{-1pt}\par\ni
In this section, we shall determine the isomorphism classes of the
Hamiltonian Lie algebras of the form $\HH=\HH(\v\ell,\G)$. We assume
that $\mbb{F}$ is an algebraically closed field.
\par
By (2.25), (2.27) and (2.35), we can rewrite (2.36) in the following more
explicit form:
$$
\begin{array}{lll}
[x^{\a,\v i},x^{\b,\v j}] &=&\dis\sum_{p\in I_{1,4}}(\a_p\b_{_{\sc
\ol p}}-\a_{_{\sc \ol p}}\b_p) x^{\si_p+\a+\b,\v i+\v j}
+\sum_{p\in I_{5,6}}(\a_pj_{_{\sc \ol p}}-i_{_{\sc \ol p}}\b_p)
x^{\si_p+\a+\b,\v i+\v j-\es_{\ol p}} \vs{4pt}\\ &&\dis+
\sum_{p\in I_{2,4}} (i_p\b_{_{\sc \ol
p}}-j_p\a_{_{\sc \ol p}}) x^{\si_p+\a+\b,\v i+\v j-\es_p} \vs{4pt}\\
&&\dis + \sum_{p\in I_4\cup I_{6,7}}(i_pj_{_{\sc \ol p}}-i_{_{\sc
\ol p}}j_p) x^{\si_p+\a+\b,\v i+\v j-\es_p-\es_{\ol p}},
\end{array}
\eqno(3.1)$$ for $(\a,\v i),(\b,\v j)\in\G\times\JJ$, where the
first summand over $p\in I_{1,4}$ corresponds to the fact that
$T_i\ne d\ne T_{\ol i}$ for $i=1,2,3,4$ (cf.~(2.12)). As for other
summands in (3.1), they are also obvious by (2.12). In particular,
we have
$$
[x^\a,x^\b]=\sum_{p\in I_{1,4}}(\a_p\b_{_{\sc \ol p}} -\a_{_{\sc
\ol p}}\b_p) x^{\si_p+\a+\b}=\sum_{p\in I_{1,4}}
\left|\matrix{\a_{\{p,\ol p\}}\cr\b_{\{p,\ol
p\}}}\right|x^{\si_p+\a+\b} \ \ \ \for\ \ \a,\b\in\G, \eqno(3.2)$$
(cf.~(2.31) and (2.22)), where $\left|\matrix{\a_{\{p,\ol
p\}}\cr\b_{\{p,\ol p\}}}\right|= \left|\matrix{\a_p&\!\!\a_{\ol
p}\cr\b_p&\!\!\b_{\ol p}}\right|$ is a $2\times2$ determinant, and
$$
[x^{-\si_p},x^{\b,\v j}]=\left\{\matrix{
(\b_p-\b_{\ol p})x^{\b,\v j}\hfill&\mbox{if \ \ \ }p\in I_1,
\vs{4pt}\hfill\cr
-\b_{\ol p}x^{\b,\v j}\hfill&\mbox{if \ \ \ }p\in I_2,
\vs{4pt}\hfill\cr
(\b_p-\b_{\ol p})x^{\b,\v j}+j_px^{\b,\v j-\es_p}
\hfill&\mbox{if \ \ \ }p\in I_3,
\vs{4pt}\hfill\cr
(\b_p-\b_{\ol p})x^{\b,\v j}+j_px^{\b,\v j-\es_p}
-j_{\ol p}x^{\b,\v j-\es_{\ol p}}
\hfill&\mbox{if \ \ \ }p\in I_4,
\hfill\cr}\right.
\eqno(3.3)$$
and
$$
[t_{\ol q},x^{\b,\v j}]=\left\{\matrix{
-\b_qx^{\b,\v j}\hfill&\mbox{if \ \ \ }q\in I_5,
\vs{4pt}\hfill\cr
-\b_qx^{\b,\v j}-j_qx^{\b,\v j-\es_q}\hfill&\mbox{if \ \ \ }q\in I_6.
\hfill\cr}\right.
\eqno(3.4)$$
\par
For any $\v i\in\JJ$, we define the {\it level} of $\v i$ to be
$$
|\v i|=\sum_{p\in J}i_p.
\eqno(3.5)$$
For any $(\a,\v i)\in\G\times\JJ$, we define the {\it support}
of $(\a,\v i)$ to be
$$
{\rm supp}(\a,\v i)=\{p\in J\,|\,\a_p\ne0\mbox{ or }i_p\ne0\}.
\eqno(3.6)$$
\par
For any Lie algebra ${\cal L}$, we denote by ${\cal L}^F$ and
by ${\cal L}^N$ the sets
of {\it ad}-locally finite elements and of {\it ad}-locally nilpotent
elements, of $\cal L$ respectively.
Generally, to obtain the isomorphism
theorem, the ordinary way is first to find the sets
$\HH^F$ and
$\HH^N$. However,
in our case here,
the determinations of the sets $\HH^F$ and $\HH^N$
seem to be un-achievable. Thus, we use a
``sandwich'' method to estimate them.
To do this, we introduce the following three subsets of $\HH$.
Denote
$$
H_1=\{x^{-\si_p},t_{\ol q}\,|\,p\in I_{1,4},q\in I_{5,6}\},
\eqno(3.7)$$
$$
H_2=\{x^{\a,\v i}\,|\,\a_{_{\sc J_{1,4}}}= \v i_{J_{1,4}\cup\ol
I_{5,6}}=0,\ i_pi_{\ol p}=0\ \for\ p\in J_7\}, \eqno(3.8)$$
$$
H_3={\rm span}\{x^{\a,\v i}\,|\,\a_{_{\sc J_{1,4}}}= \v
i_{J_{1,4}\cup\ol I_{5,6}}=0\}, \eqno(3.9)$$ (cf.~(2.22). Then our
first result is the following ``sandwich'' lemma.
\par\ni
{\bf Lemma 3.1}. {\it
$$
H_1\cup H_2\subset\HH^F\subset{\it span}(H_1\cup H_3),
\eqno(3.10)$$
$$
H_2\subset\HH^N\subset H_3.
\eqno(3.11)$$
}
{\it Proof}.
By (3.3) and (3.4), we have $H_1\subset\HH^F$.
Suppose $x^{\a,\v i}\in H_2$.
Then by (3.8),
$$
{\rm supp}(\a,\v i)\subset I_{5,6}\cup J_7,\mbox{ \ \ \ and  \ \ \
} \ol p\notin{\rm supp}(\a,\v i)\mbox{ \ if \ } p\in{\rm
supp}(\a,\v i). \eqno(3.12)$$ Let $x^{\b,\v j}\in\HH$. By (3.1)
and (3.12), we see
$$
\matrix{
[x^{\a,\v i},x^{\b,\v j}]=\!\!\!\!&
0\mbox{ \ or
a linear combination of the elements \ }x^{\g,\v k}\mbox{ \ such that}
\vs{4pt}\hfill\cr&
\mbox{there exists at least a \ }
p\in (\ol I_{5,6}\cup J_7)\bs{\rm supp}(\a,\v i)
\mbox{ \  with \ }k_p<j_p.
\hfill\cr}
\eqno(3.13)$$
Thus if we set
$$
m=1+\sum_{p\in (\ol I_{5,6}\cup J_7)\bs{\rm supp}(\a,\v i)}j_p,
\eqno(3.14)$$
then $\ad^m_{x^{\a,\v i}}(x^{\b,\v j})=0$. This proves
$H_2\subset\HH^N\subset\HH^F$.
\par
Suppose $u\notin{\rm span}(H_1\cup H_3)$.
Write
$$
u=\sum_{(\a,\v i)\in S_0}c_{\a,\v i}x^{\a,\v i},
\mbox{ \ \ where}
\eqno(3.15)$$
$$
S_0=\{(\a,\v i)\in\G\times\JJ\,|\,c_{\a,\v i}\ne0\}
\mbox{ \ \ is a finite set}.
\eqno(3.16)$$
Then by (3.7) and (3.9), there exist $(\g,\v k)\in S_0$ and
$p\in I_{1,6}$ such that at least one of $p$ and $\ol p$ is in
${\rm supp}(\g,\v k)$, mainly,
$$
(\g_p,\g_{\ol p},k_p,k_{\ol p})\ne0,
\eqno(3.17)$$
and such that
$$
(\g,\v k)\ne(-\si_p,0)\mbox{ \ \ \ if \ \ \ }p\in I_{1,4},\mbox{ \ \ \ and}
\eqno(3.18)$$
$$
(\g,\v k)\ne(0,\es_{\ol p}),\;\;k_{\ol p}\ne0
\mbox{ \ \ \ if \ \ \ }p\in I_{5,6}.
\eqno(3.19)$$
We prove that $u$ is not {\it ad}-locally finite. To do this,
we choose a total order on $\G$ compatible with group
structure of $\G$ and define the total order on $\G\times\JJ$ by
the lexicographical order,
such that the maximal element $(\g,\v k)$ of $S_0$ satisfies
(3.17)-(3.19) for some $p\in I_{1,6}$, and that
$\si_p>\si_q$ for all $q\ne p$. This is possible
because the set of all nonzero $\si_q$ is $\mbb{F}$-linear independent.
To see how it works, say, $p\in I_1$ and $(\g_p,\g_{\ol p})\ne0$
(the proof for other cases is similar).
Choose
$\b=b{\ssc\,}\es_{\ol p}\in\G$ for some $b\in\mbb{F}\bs\{0\}$ (cf.~(2.26))
such that
$$
\g_pb+m(\g_{\ol p}-\g_p)\ne0\mbox{ \ \ \  for all \ \ \ }m\in\mbb{N}.
\eqno(3.20)$$
Then for $n\in\mbb{N}$,
the ``highest'' term of $\ad_u^n(x^{\b})$ is
$x^{\b+n\g+n\si_p,n\v k}$ with the coefficient
$$
\prod_{m=0}^{n-1}(\g_p(\b_{\ol p}+m\g_{\ol p}-m)-
\g_{\ol p}(m\g_p-m))=
\prod_{m=0}^{n-1}(\g_pb+m(\g_{\ol p}-\g_p)\ne0.
\eqno(3.21)$$
Thus by (3.18), the set $\{\ad_u^n(x^{\b})\,|\,n\in\mbb{N}\}$ is linearly
independent, which implies
$$
{\rm dim}({\rm span}\{\ad_u^n(x^\b)\,|\,n\in\mbb{N}\})=\infty.
\eqno(3.22)$$
Thus $u\notin\HH^F$. This proves
$\HH^F\subset{\rm span}(H_1\cup H_3)$. Similarly,
$\HH^N\subset H_3$.
\qed\par
For any subset $X\subset\HH$, we denote by $E(X)$ the set of
the zero vector and the common eigenvectors in $\HH$ for $\ad_X$, mainly
$$
E(X)=\{u\in\HH\,|\,[X,u]\subset\mbb{F}u\}.
\eqno(3.23)$$
\par
Next, we shall determine $E(\HH^F)$. To this end, we need to
find the eigenvalues for elements of $\ad_{H_1}$.
So we define a map $\pi:\G\rar\mbb{F}^{\iota_6}$ by
$$
\pi(\a)=\mu=
(\mu_1,\cdots,\mu_{\iota_6}),
\mbox{ \ \ \ with}
\eqno(3.24)$$
$$
\mu_p=\left\{\matrix{
\a_p-\a_{\ol p}\hfill&\mbox{if }p\in I_1\cup I_{3,4},\vs{4pt}\hfill\cr
-\a_{\ol p}\hfill&\mbox{if }p\in I_2,\vs{4pt}\hfill\cr
-\a_p\hfill&\mbox{if }p\in I_{5,6},
\hfill\cr
}\right.
\eqno(3.25)$$
(cf.~(3.3) and (3.4)).
We define
$$
\MM={\rm span}\{x^\a\in\HH\,|\,\a\in\G\},
\eqno(3.26)$$
$$
\MM_\mu={\rm span}\{x^\a\,|\,\pi(\a)=\mu\}\;\;\;\;\for\;\;\;\;\mu\in\pi(\G).
\eqno(3.27)$$
Then we have
\par\ni
{\bf Lemma 3.2}. {\it
$$
E(\HH^F)=\bigcup_{\mu\in\pi(\G)}\MM_\mu,
\eqno(3.28)$$
thus $\MM={\it span}(E(\HH^F))$.
}\par\ni
{\it Proof}. By (3.10) and the definition (3.23), we have
$$
E(H_1\cup H_2)\supset E(\HH^F)\supset
E({\rm span}(H_1\cup H_3)).
\eqno(3.29)$$
We want to prove
$$
E(H_1\cup H_2)\subset
\bigcup_{\mu\in\pi(\G)}\MM_\mu
\subset E({\rm span}(H_1\cup H_3)).
\eqno(3.30)$$
Let $\mu\in\pi(\G)$. By (3.3), (3.4), (3.7)-(3.9) and (3.24)-(3.27),
elements in $\MM_\mu$ are common eigenvectors for $\ad_{H_1}$, and
$\ad_{H_3}$ acts trivially on $\MM_\mu$. Since
elements in $H_1$ commute with each other, elements in $\MM_\mu$
are common eigenvectors for $\ad_{{\rm span}(H_1\cup H_3)}$.
That is,
$$
\bigcup_{\mu\in\pi(\G)}\MM_\mu
\subset E({\rm span}(H_1\cup H_3)).
\eqno(3.31)$$
Suppose
$$
u=\sum_{(\a,\v i)\in S_0} c_{\a,\v i}x^{\a,\v i}\in\HH,
\;\;\mbox{ where }\;\;S_0=\{(\a,\v i)\in\G\times\JJ\,|\,c_{\a,\v i}\ne0\},
\eqno(3.32)$$
is a common eigenvector for $\ad_{H_1\cup H_2}$.
Since $\ad_{H_2}$ is locally nilpotent, $\ad_{H_2}$ must act trivially
on $u$. If $(\a,\v i)\in S_0$ with $i_p\ne0$ for some
$p\in\ol I_{5,6}\cup J_7$, then we can choose $v\in H_2$:
$$
v=\left\{\matrix{ x^{\es_{\ol p}}\hfill&\mbox{if}&p\in \ol
I_{5,6}, \vs{4pt}\hfill\cr t_{\ol p}\hfill&\mbox{if}&p\in J_7,
\hfill\cr}\right. \eqno(3.33)$$ such that $[v,x^{\a,\v i}]\ne0$ by
(3.1) and thus $[v,u]\ne0$, contradicting the fact that
$\ad_{H_2}$ acts trivially on $u$. Thus $\v i_{\ol I_{5,6}\cup
J_7}=0$. Similarly, since $u$ is a common eigenvector for
$\ad_{H_1}$, we must have $\v i_{I_{2,3}\cup J_4}=0$ (and thus $\v
i=0$) and $\pi(\a)=\mu$ for some $\mu$ if $(\a,\v i)\in S_0$.
This shows that $u\in\MM_\mu$. This together with (3.31) proves
(3.30). Now (3.29) and (3.30) show that all these sets are equal,
i.e., we have (3.28). \qed\par Next we shall determine the sets
$\MM^F$ and $\MM^N$. Recall that the Lie bracket in $\MM$ has the
simple form (3.2).
\par\ni
{\bf Lemma 3.3}.
{\it
$$
\MM^F={\it span}\{x^{-\si_p},x^\a\,|\,p\in I_{1,4}, \a_{_{\sc
J_{1,4}}}=0\}, \eqno(3.34)$$
$$
\MM^N={\it span}\{x^\a\,|\,\a_{_{\sc J_{1,4}}}=0\}. \eqno(3.35)$$
}
{\it Proof.}
We shall prove (3.34) as the proof (3.35) is similar.
It is straightforward to verify that by (3.2) elements in the right-hand
side of (3.34) commute with each other and they
are {\it ad}-locally finite on $\MM$.
Thus the right-hand side of (3.34) is contained in $\MM^F$.
Conversely, suppose $u\in\MM$ is not in the right-hand side of (3.34).
Then we can write $u$ as in (3.15), where now
$$
S_0=\{(\a,\v i)\in\G\times\JJ\,|\,\v i=0,\,c_{\a,\v i}\ne0\}
\mbox{ \ \ is a finite set}.
\eqno(3.36)$$
Thus we still have (3.17)-(3.19), and the same arguments after (3.19) show
that $u$ is not {\it ad}-locally finite on $\MM.$
\qed\par
Now we shall study some important features of the Lie algebra $\MM$, which
is crucial in the proof of the isomorphism theorem.
\par\ni
{\bf Lemma 3.4}. {\it (1) Assume that $\iota_4\ne0$. For
$\mu\in\pi(\G)$, regarding $\MM_\mu$ as an $\MM_0$-module via the
adjoint action, we have (i) if $\mu_{I_{1,4}}=0$, then the action of
$\MM_0$ on $\MM_\mu$ is trivial and (ii) if $\mu_{I_{1,4}}\ne0$, then
$\MM_\mu$ is a cyclic $\MM_0$-module, the nonzero multiplicative
scalars of $x^\a$ for all $\a\in\G$ with $\pi(\a)=\mu$, are the
only generators.
\par
(2) Assume that $\iota_4=0$ and $\iota_6\ne0$. Then
$(\cup_{\a\in\G}\mbb{F} x^\a)\bs\{0\}$ are the set of the common eigenvectors
of $\HH^F$ in $\MM$.
}\par\ni
{\it Proof}. (1)
Assume that $\iota_4\ne0$.
{}From (3.2) and the definition of $\pi$ in (3.24),
we see that $x^\a$ commutes with
$x^\b$ if $\pi(\a)=0$ and $(\pi(\b))_{I_{1,4}}=0$. Thus
if $\mu_{ I_{1,4}}=0$, the adjoint action of $\MM_0$ on $\MM_\mu$
is trivial.
Assume
$$
u=\sum_{\b\in S_0}c_\b x^\b\in\MM_\mu\mbox{ \ \ with \ \ }
\mu_{I_{1,4}}\ne0,\mbox{ \ \ where}
\eqno(3.37)$$
$$
S_0=\{\b\in\G\,|\,\pi(\b)=\mu,c_\b\ne0\}\mbox{ \ \ is a finite set.}
\eqno(3.38)$$
By (3.2), one has
$$
[x^\a,u]=-\sum_{p\in I_{1,4}}
\a_p\mu_p x^{\si_p+\a}\cdot u\mbox{ \ \ \ if \ \ \ }\pi(\a)=0.
\eqno(3.39)$$
Thus the subspace
$$
U={\rm span}\{x^{\si_p+\a}\cdot u= \sum_{\b\in S_0}c_\b
x^{\si_p+\a+\b}\, |\,\a\in\kn_\pi,p\in I_{1,4}\}, \eqno(3.40)$$ is
a $\MM_0$-submodule of $\MM_\mu$. Let $\la u\ra$ denote the cyclic
submodule of $\MM_\mu$ generated by $u$. Then $\la u\ra\subset U$.
If the size $|S_0|$ of $S_0$ is $\ge2$, then $U$ in (3.40) is a
proper submodule of $\MM_\mu$ and so $u$ is not a generator of
$\MM_\mu$.
\par
Now assume that $S_0$ is a singleton $\{\b\}$ with $\pi(\b)=\mu$.
Suppose $\mu_p\ne0$ for some $p\in I_{1,4}$.
For any $k\ne1$, by (3.25), $k\si_p\in\kn_\pi$, thus
$$
x^{\b+k\si_p}=-((k-1)\mu_p)^{-1}[x^{(k-1)\si_p},x^\b]\ \in\ \la u\ra.
\eqno(3.41)$$
For any $\a\in\kn_\pi$, by (3.25), $\a-(k+1)\si_p\in\kn_\pi$.
Thus by (3.2), (3.25) and (3.41), noting that $\b_{\ol q}=\b_q-\mu_q$
for $q\in I_{1,4}$,
it is straightforward to compute
that
$$
k\mu_px^{\a+\b}+\sum_{q\in I_{1,4}}(\d_{p,q}-\a_q)\mu_qx^{\a+\b-\si_p+\si_q}
=[x^{\a-(k+1)\si_p},x^{\b+k\si_p}]\ \in\ \la u\ra.
\eqno(3.42)$$
This shows that $x^{\a+\b}\in\la u\ra$ for all $\a
\in\kn_{\pi_1}$, but $\MM_\mu$ is spanned by such elements.
Thus $u$ is a generator of $\MM_\mu$.
\par
(2) is obtained directly from (3.28).
\qed\par
Let $\HH(\v\ell',\G')$ be another Hamiltonian Lie algebra defined in
last section. We shall add a prime on all the constructional
ingredients related to $\HH(\v\ell',\G')$; for instance,
$\HH',\;\JJ',\;\si'_i,\;\ell_i',\;\iota_i',$ etc.
\par
To state our isomorphism theorem, denote by $M_{m\times
n}\BF$ the space of $m\times n$ matrices with entries in
$\mbb{F}$ and by $GL_m\BF$ the group of $m\times m$
invertible matrices with entries in $\mbb{F}$.
\par\ni
{\it Definition 3.5}. Let $\G,\G'$ be two additive subgroups of
$\mbb{F}^{2\iota_7}$ satisfying (2.25) and (2.26). A group
isomorphism $\tau:\a\mapsto\a^*$ from $\G\to\G'$ is called {\it
preserving} if $\tau$ has the following form: there exists a
permutation $\nu:p\mapsto p^*$ on the index set $I_{1,4}$, which
maps $I_k\to I_k$ for $k=1,2,3,4$, such that
$$
\a^*_{\{p^*,\ol{p^*}{\ssc\,}\}}=\a_{\{p,\ol p{\ssc\,}\}}A_p\;\;\;
\for\;\;\;p\in I_{1,4}, \eqno(3.43)$$ (cf.~(2.22)), where $A_p\in
GL_2\BF$ and the multiplication in the right-hand side of (3.43)
is the vector-matrix multiplication, and
$$
A_p=\pmatrix{a_p+b_p&a_p\cr1-a_p-b_p&1-a_p\cr}\mbox{ \ or \ }
\pmatrix{1&0\cr a_p&b_p\cr}\mbox{ \ or \ }
\pmatrix{b_p&0\cr1-b_p&1\cr},
\eqno(3.44)$$
if $p\in I_1\cup I_4$
or $I_2$ or $I_3$ respectively, for some $a_p,b_p\in\mbb{F}$ with $b_p\ne0$,
and
$$
\a^*_{_{\sc I_5}}=(\a_{_{\sc I_1}}-\a_{_{\sc \ol
I_1}})B_{1,5}-\a_{_{\sc \ol I_2}} B_{2,5}+\a_{_{\sc
I_5}}B_{5,5},\mbox{ \ where} \eqno(3.45)$$
$$
B_{1,5}\in M_{\ell_1\times\ell_5}\BF,\;\;B_{2,5}\in
M_{\ell_2\times\ell_5}\BF,\;\;B_{5,5}\in GL_{\ell_5}\BF,
\eqno(3.46)$$
and
$$
\a^*_{_{\sc I_6}}=(\a_{_{\sc I_1}}-\a_{_{\sc \ol
I_1}})B_{1,6}-\a_{_{\sc \ol I_2}}B_{2,6} +(\a_{_{\sc
I_{3,4}}}-\a_{_{\sc \ol I_{3,4}}})B_{3,6}+\a_{_{\sc I_5}}B_{5,6}
+\a_{_{\sc I_6}}B_{6,6},\mbox{ \ where} \eqno(3.47)$$
$$
B_{1,6}\in M_{\ell_1\times\ell_6}\BF,\;\; B_{2,6}\in
M_{\ell_2\times\ell_6}\BF,\;\;B_{3,6}\in
M_{(\ell_3+\ell_4)\times\ell_6}\BF,\;\;
B_{5,6}\in M_{\ell_5\times\ell_6}\BF,\;\;B_{6,6}\in
GL_{\ell_6}\BF.
\eqno(3.48)$$
\par
Note that the above uniquely determine the isomorphism by (2.25).
Let us explain the above definition. First we introduce the following
notations. For any $m\times n$ matrix $A=(a_{p,q})$, we denote
by $\wt A=(\wt a_{p,q})$
(resp.~$\wh A=(\wh a_{p,q}){\sc\,}$)
the $2m\times n$ matrix such that the odd rows of $\wt A$ (resp.~$\wh A$)
forms the matrix $A$ (resp.~the $m\times n$ zero matrix)
and the even rows of $\wt A$ (resp.~$\wh A$) forms the matrix $-A$, i.e.,
$$
\wt a_{2p-1,q}=-\wt a_{2p,q}=a_{p,q},\;\;\;
\wh a_{2p-1,q}=0,\;\wh a_{2p,q}=-a_{p,q}\;\;\;
\for\;\;\;p\in\ol{1,m},\,q\in\ol{1,n}.
\eqno(3.49)$$
A preserving isomorphism $\tau$ can be
decomposed into the composition of two isomorphisms $\tau=\tau_\nu\cdot
\tau_0$ such that $\tau_\nu$ only involves the permutation $\nu$, i.e.,
in (3.43)-(3.48), all $A_p$ and $B_{i,i}$ are identity matrices and
all $B_{i,j}$ are zero matrices for $i\ne j$; and $\tau_0$ only involves
matrices, i.e., $\nu={\bf 1}_{I_{1,4}}$ in (3.43). Furthermore, $\tau_0$
can be decomposed into $\tau_0=\tau_1\cdot\tau_2$ such that $\tau_1,\tau_2$
have the following forms
$$
\tau_1:\ \ (\a^*_{_{\sc J_{1,4}}},\a^*_{_{\sc I_{5,6}}})
=(\a_{_{\sc J_{1,4}}},\a_{_{\sc I_{5,6}}})A,\ \ \mbox{ where}
\eqno(3.50)$$
$$
A={\rm diag}(A_1,\cdots,A_{\iota_4},B_{5,5},B_{6,6}),
\eqno(3.51)$$
and
$$
\tau_2:\ \ (\a^*_{_{\sc J_{1,4}}},\a^*_{_{\sc I_{5,6}}})
=(\a_{_{\sc J_{1,4}}},\a_{_{\sc I_{5,6}}})C,\ \ \ C={\bf
1}_{2\iota_4+\ell_5+\ell_6}+D, \eqno(3.52)$$ where in general
${\bf1}_m$ denotes the $m\times m$ identity matrix, and where $D$
has the form
$$
D=(0,D_5,D_6),\ \
D_5=
\left(\begin{array}{l}\wt B_{1,5}\\ \wh B_{2,5}\\
0\end{array}\right),\ \
D_6=
\left(\begin{array}{l}\wt B_{1,6}\\ \wh B_{2,6}\\ B_{3,6}\\ B_{5,6}\\
0\end{array}\right),
\eqno(3.53)$$
where $0$ denotes some proper zero matrices whose orders are
clear from the context.
\par
Now we can state the main result of this paper.
\par\ni
{\bf Theorem 3.6}. {\it $\th:\HH(\v\ell,\G)\cong\HH(\v \ell',\G')$ if and only
if $\v\ell=\v\ell'$ and there exists a preserving isomorphism
$\tau:\G\cong\G'$.}
\par\ni
{\bf Theorem 3.7 (Main Theorem)}. {\it Two Hamiltonian Lie
algebras are isomorphic if and only if their corresponding Poisson
algebras are isomorphic. }\par\ni {\it Proof}. By Theorem 3.6 and
by [SX], the condition for two Hamiltonian Lie algebras being
isomorphic is the same as the condition for the corresponding two
Poisson algebras being isomorphic. \qed\par\ni {\it Proof of
Theorem 3.6}. ``$\Leftarrow$'': Suppose $\v\ell=\v\ell'$ and
$\tau:\G\to\G'$ is a preserving isomorphism. By the explanation
above, $\tau$ can be written as
$\tau=\tau_\nu\cdot\tau_1\cdot\tau_2$, thus it suffices to
consider the following 3 cases.
\par
{\it Case a}: First assume that $\tau=\tau_\nu$ is determined by
permutation $\nu$.
\par
For any $\v i\in\JJ$, we define $\v i^*\in\JJ$ which is obtained from
$\v i$ by permutation $\nu$. Then it is straightforward to verify that
the linear map
$$
\th_\nu:\HH\to\HH'\mbox{ \ such that \ }\th_\nu(x^{\a,\v
i})=x^{\a^*,\v i^*}, \eqno(3.54)$$ is a Lie algebra isomorphism.
\par
{\it Case b}: Next assume that $\tau=\tau_1$ as in (3.50).
\par
We shall define an isomorphism $\th:\HH\to\HH'$ as Poisson algebra
isomorphism (then $\th$ is clearly a Lie algebra isomorphism).
By (1.1), it suffices to find the images of
the generators $x^\a,t_p$ for $\a\in\G,p\in I_{2,4}\cup I_{6,7}
\cup\ol I_{4,7}$ (cf.~(3.58) and (3.62)-(3.64) below)
such that the following conditions hold
(cf.~[SX]):
$$
\th([x^\a,x^\b])=[\th(x^\a),\th(x^\b)],\;\;\;
\th([t_p,x^\b])=[\th(t_p),\th(x^\b)],\;\;\;
\th([t_p,t_q])=[\th(t_p),\th(t_q)],
\eqno(3.55)$$
for $\a,\b\in\G$ and $p,q\in I_{2,4}\cup I_{6,7}\cup\ol I_{4,7}$.
\par
Let $\D=\sum_{p\in I_{1,4}}\mb{Z}\si_p$ be the subgroup of $\G$ generated by
$\{\si_p\,|\,p\in I_{1,4}\}$ and define $\chi:\D\to\mbb{F}^\times=
\mbb{F}\bs\{0\}$ to be
the {\it character} of $\D$ (i.e., the group homomorphism
$\D\to\mbb{F}^\times$) determined by
$$
\chi(\si_p)=b_p\;\;\;\for\;\;\;p\in I_{1,4},
\eqno(3.56)$$
where $b_p$ are elements in $\mbb{F}$ appearing as entries of matrices $A_p$
in (3.44). We prove that $\chi$ can be extended to a character
$\chi:\G\to\mbb{F}^\times$ as follows: Assume that $\D_1\supset\D$ is a
maximal subgroup of $\G$ such that $\chi$ can be extended to a character
$\chi:\D_1\to\mbb{F}^\times$. If $\D_1\ne\G$, then we choose $\a\in\G\bs\D_1$
and extend $\chi$ to $\D_2=\mb{Z}\a+\D_1\to\mbb{F}^\times$ by defining
$$
\chi(m\a+\b)=\left\{\begin{array}{ll}
\chi(\b)&\mbox{if \ }\mb{Z}\a\cap\D_1=\{0\},
\vs{4pt}\\
a^m\chi(\b)&\mbox{if \ }\mb{Z}\a\cap\D_1=\mb{Z}n\a,
\end{array}\right.
\eqno(3.57)$$
for $m\in\mb{Z},\b\in\D_1$, where $a$ is an $n$th root of $\chi(n\a)$ in the
second case (recall that $\mbb{F}$ is algebraically closed). This leads to
a contradiction with the maximality of $\D_1$. Thus
$\chi$ can be extended to a character $\chi:\G\to\mbb{F}^\times$.
\par
Now we define the images of $x^\a$ to be
$$
\th(x^\a)=\chi(\a)x'{}^{\a^*}\;\;\;\;\for\;\;\;\;\a\in\G,
\eqno(3.58)$$
(recall that we add prime on the constructional ingredients related to
$\HH'$). Then by (3.2) we see that the first equation of (3.55) holds because
(3.44) and
(3.50) guarantees that $\si_p^*=\si'_p$ and that the determinant of $A_p$
is $|A_p|=b_p=\chi(\si_p)$ and
$$
\chi(\a)\chi(\b) \left|\matrix{\a^*_{\{p,\ol p\}}\cr\b^*_{\{p,\ol
p\}}}\right| =\chi(\a+\b) \left|\matrix{\a_{\{p,\ol
p\}}\cr\b_{\{p,\ol p\}}}\right|\cdot|A_p|
=\chi(\si_p+\a+\b)\left|\matrix{\a_{\{p,\ol p\}}\cr\b_{\{p,\ol
p\}}}\right|. \eqno(3.59)$$
 Next we shall find the image of $t_p$.
To do this, we introduce a new notation: For any vector
$s=(s_1,s_{\ol1},s_2,s_{\ol2},\cdots, s_{\iota_7},s_{\ol\iota_7})$
(with entries in $\mbb{F},\HH$ or in $\HH'$), we denote
$$
\ol s=(-s_{\ol1},s_1,-s_{\ol2},s_2,\cdots,-s_{\ol\iota_7},s_{\iota_7}).
\eqno(3.60)$$
For a subset $K\subset J$, we denote by $\ol s_K$ the vector obtained from
$\ol s$ by deleting $-s_{\ol p},s_q$ for $\ol p,q\in J\bs K$; for instance,
$$
\ol s_{\{\ol1,\ol2,3,4,\ol4\}}=(-s_{\ol1},-s_{\ol2},s_3,-s_{\ol4},s_4),
\eqno(3.61)$$
(cf.~(2.22)).
We define
$$
\th(t_p)=s_p\;\;\;\for\;\;\;p\in I_{2,4}\cup I_{6,7}\cup\ol I_{4,7},
\mbox{ \ where}
\eqno(3.62)$$
$$
s_p=t'_p,\;\;\; s_q=b_qt'_q,\;\;\;
 (-s_{\ol r},s_r)=b_r(-t'_{\ol
r},t'_r)A_r^{-1}\;\;\; \for\;\;\;p\in I_2,\;q\in I_3,\;r\in I_4,
\eqno(3.63)$$
$$
\ol s_{_{\sc \ol I_{5,6}}}=\ol t'_{_{\sc \ol I_{5,6}}} {\rm
diag}(B_{5,5},B_{6,6})^{-1},\;\;\; \ol s_{_{\sc I_6}}=\ol
t'_{_{\sc I_6}}B_{6,6}^T,\;\;\;
s_{_{\sc J_7}}=t'_{_{\sc J_7}},
\eqno(3.64)$$ where the up-index ``T'' stands for the transpose of
a matrix. Then if $p\in I_{2,3}$, we have
$$
[\th(t_p),\th(x^\a)]=\chi(\a)b_p \a_{\ol
p}x'{}^{\a^*+\si'_p}=\th(\a_{\ol p}x^{\a+\si_p}) =\th([t_p,x^\a]),
\eqno(3.65)$$ because by (3.44) and (3.50),
 $\a^*_{\ol
p}=b_p\a_{\ol p}$ if $p\in I_2$ and $\a^*_{\ol p}=\a_{\ol p}$ if
$p\in I_3$. If $p\in I_4$, as $1\times2$ matrices with entries in
$\HH$, we have
$$
[\th(\ol t_{\{p,\ol p\}}),\th(x^\a)]= \chi(\a)b_p\a^*_{\{p,\ol
p\}}A_p^{-1}x'{}^{\a^*+\si'_p}= \chi(\a+\si_p)\a_{\{p,\ol
p\}}x'{}^{\a^*+\si'_p}= \th([\ol t_{\{p,\ol p\}},x^\a]).
\eqno(3.66)$$ Furthermore, we have $[(\ol t_{_{\sc \ol
I_{5,6}}},\ol t_{_{\sc I_6\cup J_7}}),x^\a]= (\a_{_{\sc
I_{5,6}}},0)x^\a$, and
$$
\a^*_{_{\sc I_{5,6}}}=\a_{_{\sc I_{5,6}}}{\rm
diag}(B_{5,5},B_{6,6}), \eqno(3.67)$$ by (3.50). {}From this and
(3.64), we obtain
$$
[(\th(\ol t_{_{\sc \ol I_{5,6}}}), \th(\ol t_{_{\sc I_6\cup
J_7}})),\th(x^\a)] =\th([(\ol t_{_{\sc \ol I_{5,6}}},\ol t_{_{\sc
I_6\cup J_7}}),x^\a]). \eqno(3.68)$$ {}From this and (3.66), we
obtain the second equation of (3.55).
\par
To verify the last equation of (3.55), note that
$$
[\ol t_{_{\sc I_{2,3}\cup J_4\cup I_5\cup J_{6,7}}}^T, \ol
t_{_{\sc I_{2,3}\cup J_4\cup I_5\cup J_{6,7}}}] ={\rm
diag}(0,S^\si_{\ell_4},0,S_{\ell_6+\ell_7}), \eqno(3.69)$$ where
$$
S^\si_{\ell_4}={\rm
diag}\left(\left(\matrix{0&x^{\si_{\iota_3+1}}\cr
-x^{\si_{\iota_3+1}}&0}\right),\cdots,
\left(\matrix{0&x^{\si_{\iota_4}}\cr-x^{\si_{\iota_4}}&0}
\right)\right), \eqno(3.70)$$ is a $2\ell_4\times2\ell_4$ matrix
with entries in $\HH$, and where, in general
$$
S_m={\rm diag}\left(\left(\matrix{0&1\cr
-1&0}\right),\cdots,\left(\matrix{0&1\cr-1&0}\right)\right)\in GL_{2m}.
\eqno(3.71)$$
Using (3.69), (3.63) and (3.64), we can obtain
$$
[\th(\ol t_{_{\sc I_{2,3}\cup J_4\cup I_5\cup J_{6,7}}})^T,
\th(\ol t_{_{\sc I_{2,3}\cup J_4\cup I_5\cup J_{6,7}}})] =
\th([\ol t_{_{\sc I_{2,3}\cup J_4\cup I_5\cup J_{6,7}}}^T, \ol
t_{_{\sc I_{2,3}\cup J_4\cup I_5\cup J_{6,7}}}]). \eqno(3.72)$$
For example, if $p\in I_4$, by (3.56), (3.58) and (3.63), we have
$$
\matrix{
[\th(\ol t_{\{p,\ol p\}})^T,\th(\ol t_{\{p,\ol p\}})]
\!\!\!\!&
=b_p(A_p^{-1})^T[\ol t{}^{{\sc\,}'T}_{\{p,\ol p\}},\ol t'_{\{p,\ol p\}}]b_pA_p^{-1}
\vs{4pt}\hfill\cr&
=b_p\pmatrix{0&x'{}^{\si'_p}\cr-x'{}^{\si'_p}&0\cr}
=\th([\ol t^T_{\{p,\ol p\}},\ol t_{\{p,\ol p\}}]).
\hfill\cr}
\eqno(3.73)$$
This proves the last equation of (3.55).
\par
{\it Case c}: Assume that $\tau=\tau_2$ as in (3.52).
\par
We define (3.58) with $\chi(\a)=1$ and
we define (3.62) with
$$
\ol s_{_{\sc I_{2,3}\cup J_4}}=\ol t'_{_{\sc I_{2,3}\cup J_4}}+\ol
t'_{_{\sc I_6}}E_1, \eqno(3.74)$$
$$
\ol s_{_{\sc \ol I_{5,6}}}= \ol t'_{_{\sc \ol I_{5,6}}}E_2+\ol
t'_{_{\sc I_6}}E_3+x'{}^{-\si'}E_4, \;\;\; \ol s_{_{\sc I_6\cup
J_7}}=\ol t'_{_{\sc I_6\cup J_7}}, \eqno(3.75)$$ where
$E_1,...,E_4$ are some matrices to be determined in order that
(3.55) holds and where $x'{}^{-\si'}$ denotes the vector
$$
x'{}^{-\si'}=(x'{}^{-\si'_1},\cdots,x'{}^{-\si'_{\iota_4}}).
\eqno(3.76)$$ We shall not give the explicit forms of
$E_1,...,E_4$ here, but an interested reader can find the
solutions by considering two special cases of (3.53): (1) $D_5=0$,
(2) $D_6=0$ (the general case is the composition of the two
special cases), or refer to [SX] (also, cf.~the proof of
necessity).
\par
``$\Rightarrow$'': Assume that there exits a Hamiltonian
Lie algebra isomorphism $\th: \HH(\v\ell,\G)\rar\HH(\v \ell',\G')$.
\par
First, we make the following conventions: If a subset of $\HH$ is
defined, then we take the definition of the corresponding subset
of $\HH'$ for granted. If a property about $\HH$ is given, the
same property also holds for $\HH'$, without description.
\par
Clearly, $\th$ maps $\HH^F,\HH^N$ to $\HH'^F,\HH'^N$
respectively, thus also maps $\MM\to\MM'$ by Lemma 3.2. By Lemma 3.3,
we have ${\rm dim}(\MM^F/\MM^N)=\iota_4$. This shows that
$$
\iota_4=\iota'_4.
\eqno(3.77)$$
For simplicity, we assume that $\iota_4\ne0$ (if $\iota_4=0$, using Lemma
3.4 (2), one sees that all statements or arguments below either work or
do not apply to the case; if $\iota_6=0$, then one can go directly to
Claim 8 below). Denote
$$
\G_{1,4}=\{\a\in\G\,|\,(\pi(\a))_{I_{1,4}}=0\},
\eqno(3.78)$$
(cf.~(3.24) and (2.22)). By Lemma 3.2, there exists a bijection
$\tau_1:\pi(\G)\to\pi(\G')$ such that
$$
\th(\MM_\mu)=\MM'_{\tau_1(\mu)}\;\;\;\;\for\;\;\;\;\mu\in\pi(\G)\mbox{,
\ \ \ and \ }\tau_1(0)=0.
 \eqno(3.79)$$ {}From this and Lemma 3.4, there
exists a bijection $\G\bs\G_{1,4}\to \G'\bs\G'_{1,4}$ which shall
be denoted by $\tau:\a\mapsto\a^*$ such that
$$
\th(x^\a)=c_\a x'{}^{\a^*}\;\;\;\;\for\;\;\;\a\in\G\bs\G_{1,4}
\mbox{ \ \ and some \ }c_\a\in\mbb{F}^\times. \eqno(3.80)$$ We
shall prove the necessity by establishing several claims.
\par
{\bf Claim 1}.
There exists a bijection $I_{1,4}\rar I'_{1,4}$ denoted by $\nu:p\mapsto p^*$
such that
$$
\th(x^{-\si_p})=d_px'{}^{-\si'_{p^*}}\mbox{ for }p\in I_{1,4}
\mbox{ and some }d_p\in\mbb{F}^\times.
\eqno(3.81)$$
\par
By (3.7)-(3.9) and Lemma 3.3, we have
$$
\{u\in\MM^F\,|\,[u,H_1\cup H_2]=0\}=
{\rm span}\{x^{-\si_p}\,|\,p\in I_{1,4}\}=
\{u\in\BB^F\,|\,[u,H_1\cup H_3]=0\}.
\eqno(3.82)$$
Thus by Lemma 3.1,
$$
\{u\in\MM^F\,|\,[u,\HH^F]=0\}=
{\rm span}\{x^{-\si_p}\,|\,p\in I_{1,4}\}.
\eqno(3.83)$$
Let $p\in I_{1,4}$. Then by (3.83),
$$
\th(x^{-\si_p})\in\sum_{q\in I'_{1,4}}\mbb{F}x'{}^{-\si'_q}.
\eqno(3.84)$$
Suppose
$$
\th(x^{-\si_p})\notin\bigcup_{q\in I'_{1,4}}\mbb{F}x'{}^{-\si'_q}.
\eqno(3.85)$$
By (2.26), there exists $a\in\mbb{F}^\times$ such that $a\es_{\ol p}\in\G$. By
(3.2), we have
$$
[x^{a\es_{\ol p}-\si_p},x^{-a\es_{\ol p}-\si_p}]=2ax^{-\si_p}.
\eqno(3.86)$$
Note that $a\es_{\ol p}-\si_p,-a\es_{\ol p}-\si_p\notin\G_{1,4}$, by (3.81),
$$
\th(x^{a\es_{\ol p}-\si_p})\in\mbb{F} x'{}^\a\bs\{0\},\;\;\;
\th(x^{-a\es_{\ol p}-\si_p})\in\mbb{F} x'{}^\b\bs\{0\}\mbox{ \ for some \ }
\a,\b\in\G'\bs\G'_{1,4}.
\eqno(3.87)$$
By (3.2), we have
$$
[x'{}^\a,x'{}^\b]=\sum_{q\in I'_{1,4}}(\a_q\b_{\ol q}-\a_{\ol q}\b_q)x'{}
^{\si'_q+\a+\b}.
\eqno(3.88)$$
By (3.84)-(3.86) and (3.88), there exist $q,r\in I'_{1,4}$ with $q\ne r$
such that $\si'_q+\a+\b=-\si'_r$. Thus
$$
\b=-\a-\si'_q-\si'_r,
\eqno(3.89)$$
and (3.88) becomes
$$
[x'{}^\a,x'{}^\b]=
(\a_q\eta'_{\ol q}+\a_{\ol q})x'{}^{-\si'_r}+
(\a_r\eta'_{\ol r}+\a_{\ol r})x'{}^{-\si'_q},
\eqno(3.90)$$
where in general, for $q\in J_{1,4}$, we denote
$$
\eta_q=\left\{\matrix{1\hfill&\mbox{if \ }q\in I_{1,4},\vs{4pt}\hfill\cr
-1\hfill&\mbox{if \ }q\in\ol I_1\cup I_{3,4},\vs{4pt}\hfill\cr
0\hfill&\mbox{if \ }q\in\ol I_2,\hfill\cr}\right.
\eqno(3.91)$$
and we define $\eta'_q$ similarly
(then $\si'_q=\es'_q-\eta'_{\ol q}\es'_{\ol q}$, cf.~(2.23)).
By (3.85), both coefficients in (3.90) are nonzero.
Since $2a\es_{\ol p}-\si_p\in\G\bs\G_{1,4}$, we have
$$
\th(x^{2a\es_{\ol p}-\si_p})\in\mbb{F}x'{}^\g\bs\{0\}\mbox{ \ for some \ }
\g\in\G'\bs\G'_{1,4}.
\eqno(3.92)$$
{}From $[x^{2a\es_{\ol p}-\si_p},x^{-a\es_{\ol p}-\si_p}]\in
\mbb{F}x^{a\es_{\ol p}-\si_p},$ it follows from (3.87)
that
$$
[x'{}^\g,x'{}^\b]\in\mbb{F}x'{}^\a.
\eqno(3.93)$$
Thus there exists $q'\in I'_{1,4}$
such that
$$
\g_{q'}\b_{\ol q'}-\g_{\ol q'}\b_{q'}\ne0\mbox{ \ and \ }
\si'_{q'}+\g+\b=\a.
\eqno(3.94)$$
Hence
$$
\g=\a-\b-\si'_{q'}=2\a+\si'_q+\si'_r-\si'_{q'},
\eqno(3.95)$$
by (3.89).
If $q\ne q'\ne r$, we deduce from (3.89) and (3.95) that
$$
\matrix{
[x'{}^\g,x'{}^\b]=\!\!\!\!&
(\a_q\eta'_{\ol q}+\a_{\ol q})x'{}^{\si'_q+\g+\b}
+(\a_r\eta'_{\ol r}+\a_{\ol r})x'{}^{\si'_r+\g+\b}
\vs{4pt}\hfill\cr&
+(\g_{q'}\b_{\ol q'}-\g_{\ol q'}\b_{q'})x'{}^{\si'_{q'}+\g+\b}
\notin\mbb{F}x'{}^\a,
\hfill\cr}
\eqno(3.96)$$
a contradiction with (3.93). Similarly, if $q'=q$ or $q'=r$, we can still
deduce a contradiction from (3.89), (3.93) and (3.95). This proves
the claim.
\par
We extend $\nu$ to $\nu:J_{1,4}\to J'_{1,4}$ such that
$\nu(\ol p)=\ol p^*$ for $p\in I_{1,4}$.
For $p\in I_{1,4}$, by (2.26),
we fix $e_p\in\mbb{F}^\times$ such that
$$
\l_p=e_p\es_{\ol p}\in\G\bs\{0\}. \eqno(3.97)$$ Then
$\l_p\notin\G_{1,4}$. Denote $\l^*_p=\tau(\l_p)$ (cf.~(3.80)).
Write
$$
\l^*_p=(\l^*_{p,1},\l^*_{p,\ol 1},\cdots,\l^*_{p,\iota'_7},
\l^*_{p,\ol\iota'_7})\in\G'\subset\mbb{F}^{2\iota'_7},
\eqno(3.98)$$ (cf.~(2.20)). For $p,q\in I_{1,4}$, applying $\th$
to $[x^{\l_p},x^{-\si_q}]=\d_{p,q}e_px^{\l_p}$, by (3.80) and
(3.81), we obtain
$$
d_q(\eta'_{\ol q^*}\l^*_{p,q^*}+\l^*_{p,\ol q^*})= \d_{p,q}e_p
\qquad\mbox{for \ }p,q\in I_{1,4}. \eqno(3.99)$$ Let $p\ne q$.
Applying $\th$ to $0=[x^{\l_p},x^{\l_q}]$ and using (3.99), we
obtain
$$
0=\l^*_{p,q^*}\l^*_{q,\ol q^*}-\l^*_{p,\ol q^*}\l^*_{q,q^*}
=\l^*_{p,q^*}(\l^*_{q,\ol q^*}+\eta'_{\ol q^*}\l^*_{q,q^*})
=\l^*_{p,q^*}d_q^{-1}e_q. \eqno(3.100)$$ The above two equations
imply
$$
\l^*_{p,q^*}=0\qquad \for\ \ p\in I_{1,4},\,q\in J_{1,4},\, q\ne
p,\ol p. \eqno(3.101)$$ Denote
$$
\G_p=(\mbb{F}\es_p+\mbb{F}\es_{\ol p})\cap\G. \eqno(3.102)$$
Exactly to the proof of (3.101), we have
$$
\a^*_{q^*}=0\;\;\;\for\;\;\;\a\in\G_p\bs\G_{1,4}, \,p,q\in
J_{1,4},q\ne p,\ol p. \eqno(3.103)$$
\par
{\bf Claim 2}.
$\tau:\a\mapsto\a^*$ can be uniquely extended to a group isomorphism
$\tau:\G\rar\G'$ such that $\si^*_p=\si'_{p^*}$ for $p\in I_{1,4}$.
\par
Noting that by (3.24), (3.25) and (3.78), $\a\notin\G_{1,4}$
implies $\a+k\si_1\notin\G_{1,4}$ for $k\in\mb{Z}$.
For any
$\a\in\G,\b\in\G_1$ with $\a,\b,\a+\b\notin\G_{1,4}$, we have
(recall (3.91))
$$
\matrix{ (\a_1(\b_{\ol1}+\eta_{\ol1})
\!\!\!\!&-\a_{\ol1}(\b_1-1))c_{\a+\b-\si_1}x'{}^{(\a+\b)^*}
\vs{4pt}\hfill\cr& =c_\a c_{\b-\si_1}(\a^*_{1^*}(\b-\si_1)^*_{\ol
1^*} -\a^*_{\ol 1^*}(\b-\si_1)^*_{1^*})
x'{}^{\si'_{1^*}+\a^*+(\b-\si_1)^*}, \hfill\cr} \eqno(3.104)$$ by
applying $\th$ to (3.2) and by (3.103). By comparing the power of
$x'$, this implies
$$
(\a+\b)^*=\si'_{1^*}+\a^*+(\b-\si_1)^* \eqno(3.105)$$ if $\a,\b$
satisfy
$$
\b\in\G_1,\;\;\a,\b,\a+\b\in\G\bs\G_{1,4}, \mbox{ \ and \
}\a_1(\b_{\ol1}+\eta_{\ol1})-\a_{\ol1}(\b_1-1)\ne0. \eqno(3.106)$$
Let $\a\in\G\bs\G_{1,4}$. We prove by induction on $|k|$ that
$$
(k\a)^*-k\a^*\in\wt\G'_1,\mbox{ where
}\wt\G'_1=\{\b\in\G'\,|\,\b_q=0 \mbox{ for }q\in J_{1,4},\,q\ne
1^*,\ol 1^*\}. \eqno(3.107)$$ Let $\g\in\G$ such that
$\g,\a+\g\notin\G_{1,4}$. We have
$$
\sum_{p\in I_{1,4}}(\a_p\g_{\ol p}-\a_{\ol p}\g_p)c_{\si_p+\a+\g}
x'{}^{(\si_p+\a+\g)^*} =c_\a c_\g\sum_{p\in I'_{1,4}}
(\a^*_{p^*}\g^*_{\ol p^*}-\a^*_{\ol p^*}\g^*_{p^*})
x'{}^{\si'_{p^*}+\a^*+\g^*}.
 \eqno(3.108)$$
  We inductively assume
that (3.107) holds for $k$ (for instance, $k=1$). Let $\g=k\a+\b$
for some suitable $\b\in\G_1$ such that condition (3.106) holds
for all the involved pairs for which we need to make use of
(3.105) in the following proof (when $\a,k$ are fixed, by (2.26),
such $\b$ exists), by (3.107) (note that we assume (3.107) holds
for $k$), (3.105) and (3.103), we see that all terms in (3.108)
vanish except the terms corresponding to $p=1$ in both sides. Thus
we obtain
$$
\matrix{ \si'_{1^*}+((k+1)\a)^*+\b^*\!\!\!\!&=
(\si_1+(k+1)\a+\b)^* \vs{4pt}\hfill\cr&
=\si'_{1^*}+\a^*+(k\a+\b)^*=
2\si'_{1^*}+\a^*+(k\a)^*+(\b-\si_1)^*, \hfill\cr} \eqno(3.109)$$
where the first and last equalities follow from (3.105) and the
second follows from (3.108). {}From this we see that (3.107) holds
for $k+1$. This proves (3.107). Now replacing $\a$ by $j\a$ (with
$j\ne0$) and $\b$ by $k\a+\b-\si_1$ in (3.108) (with suitable
$\b\in\G_1$), since (3.107) holds, we have again that all terms in
(3.108) vanish except the terms corresponding to $p=1$ in both
sides. Thus we have similar formula as in (3.109):
$$
((j+k)\a+\b)^*=2\si'_{1^*}+(j\a)^*+(k\a)^*+(\b-2\si_1)^*.
\eqno(3.110)$$ {}From this we obtain
$$
(j\a)^*+(k\a)^*=(j'\a)^*+(k'\a)^*\mbox{ if
}j+k=j'+k',\,j,k,j',k'\ne0. \eqno(3.111)$$ {}From this we obtain
$$
(j\a)^*=j\a^*\mbox{ \ \ for \ \
}\a\in\G\bs\G_{1,4},j\in\mb{Z}\bs\{0\}. \eqno(3.112)$$ For some
suitable $\b\in\G_1$, by (3.105), (3.110) and (3.112), we have
$$
\matrix{ \si'_{1^*}+(j\a+\si_1)^*+(\b-2\si_1)^*
\!\!\!\!&=((j\a+\si_1)+(\b-\si_1))^* \vs{4pt}\hfill\cr&
=(j\a+\b)^* =2\si'_{1^*}+j\a^*+(\b-2\si_1)^*. \hfill\cr}
\eqno(3.113)$$ {}From this we obtain
$$
(j\a+\si_1)^*=j\a^*+\si'_{1^*} \mbox{ \ \ for \ \
}\a\in\G\bs\G_{1,4},j\in\mb{Z}\bs\{0\}. \eqno(3.114)$$ Now take
any $\a,\g\in\G$ such that
$$
\a,\g,\a+\g\in\G\bs\G_{1,4}\mbox{ and }
\a_1\g_{\ol1}-\a_{\ol1}\g_1\ne0. \eqno(3.115)$$ Using (3.114) in
(3.108), by comparing the term $x'{}^{(\si_1+\a+\g)^*}$ in both
sides, we obtain
$$
\begin{array}{ll}\dis
(\a+\g)^*=\a^*+\g^*+\sum_{p\in
I_{1,4}}k^{(p)}_{\a,\g}(\si'_{p^*}-\si'_{1^*}), \mbox{ where}
\vs{4pt}\\ \dis k^{(p)}_{\a,\g}=0,1\mbox{ such that } \sum_{p\in
I_{1,4}}k^{(p)}_{\a,\g}\le1.
\end{array}
\eqno(3.116)$$ We claim that $(\a+\g)^*=\a^*+\g^*$
if the pairs $(\a,\g),(2\a,2\g)$ satisfy (3.115). Assume that
$k_{\a,\g}^{(q)}=1$ for some $q\in I_{1,4}$. Then we obtain
$$
\begin{array}{ll}\dis
(2\a)^*+(2\g)^*+ \sum_{p\in
I_{1,4}}k_{2\a,2\g}^{(p)}(\si'_{p^*}
-\si'_{1^*})
\!\!\!\!&=(2\a+2\g)^*=(2(\a+\g))^* =2(\a+\g)^*
\\ &\dis=2(\a^*+\g^*+\sum_{p\in I_{1,4}}k_{\a,\g}^{(p)}(\si'_{p^*}
-\si'_{1^*})),
\end{array}
\eqno(3.117)$$ from this we obtain
$k_{2\a,2\g}^{(q)}=2k_{\a,\g}^{(q)}>1$, which is a contradiction
to (3.116).
\par
For any $\a,\b,\a+\b\in\G\bs\G_{1,4}$, we can always choose
$\g\in\G\bs\G_{1,4}$ such that the pairs
$$
(\a+\b,\g),\;(2\a+2\b,2\g),\;(\a,\b+\g),\;(2\a,2\b+2\g),\;
(\b,\g),\;(2\b,2\g), \eqno(3.118)$$ satisfy (3.115). Hence
$$
(\a+\b)^*+\g^*=(\a+\b+\g)^*=\a^*+(\b+\g)^*= \a^*+\b^*+\g^*,
\eqno(3.119)$$ which  shows
$$
(\a+\b)^*=\a^*+\b^*\;\;\;\for\;\;\;\a,\,\b,\,\a+\b\in\G\bs\G_{1,4}.
\eqno(3.120)$$
 This shows that $\tau$ can be uniquely extended to
a group isomorphism $\tau:\G\rar\G'$ such that $\si^*_1=
\si'_{1^*}$ and so similarly $\si^*_p=\si'_{p^*}$ for $p\in
I_{1,4}$. This proves the claim.
\par
{\bf Claim 3}. We have $\nu(I_i)=I'_i$ for $i=1,2,3,4$. In
particular,
$(\ell_1,\ell_2,\ell_3,\ell_4)=(\ell'_1,\ell'_2,\ell'_3,\ell'_4)$,
$I_i=I'_i$ for $i=1,2,3,4$, and $\si_p=\si'_p,\,\eta_p=\eta'_p$
for $p\in I_{1,4}$ (cf.~(2.23) and (3.91)).
\par
Note that $\ad_{x^{-\si_p}}$ is a semi-simple operator on $\HH$
if and only if $p\in I_{1,2}$ (cf.~(3.3)). Thus
$$
\nu(I_{1,2})=I'_{1,2},\;\;\;\mbox{ and so \
}\nu(I_{3,4})=I'_{3,4}. \eqno(3.121)$$ Denote
$$
\matrix{ \NN\!\!\!\!&=\{u\in\HH\,|\,[u,\MM]\subset\MM\}
\vs{4pt}\hfill\cr& =\MM+{\rm span}\{x^{\a,\v i}\,|\,\a=\a_{_{\sc
I_{5,6}}},\, |\v i|=1\mbox{ or }\v i=\v i_{_{\sc I_6\cup J_7}}\},
\hfill\cr} \eqno(3.122)$$
$$
\matrix{ \NN_0\!\!\!\!&=\MM+\{u\in\NN\,|\,[x^{-\si_p},u]=0\ \for\
p\in I_{1,4}\} \vs{4pt}\hfill\cr& =\MM+{\rm
span}\{x^{\a,\es_q},x^{\a,\v j}\,|\,\a=\a_{_{\sc
I_{5,6}}},\,q\in\ol I_{5,6},\, \v j=\v j_{_{\sc I_6\cup J_7}}\},
\hfill\cr} \eqno(3.123)$$
$$
\NN_p=\MM+{\rm span}\{u\in\NN\,|\,[x^{-\si_p},u]=0\}\mbox{ for
}p\in I_{1,4}. \eqno(3.124)$$ Then $\NN_0$ is a Lie algebra and
$\NN$ is an $\NN_0$-module such that $\NN_p$ is a submodule for
$p\in I_{1,4}$. Note that the quotient module $\NN/\NN_p$ is zero
if $p\in I_1$, is a cyclic $\NN_0$-module (with generator $t_{\ol
p}$) if $p\in I_{2,3}$, and is not cyclic (with two generators
$t_p,t_{\ol p}$) if $p\in I_4$. Applying $\th$ to the above sets
and by (3.121), we obtain the claim.
\par
Using Claim 3 and (3.54), by replacing $\HH$ by $\th_\nu(\HH)$
(cf.~(3.54)), we can now suppose $\nu=1$.
\par
{\bf Claim 4}.
There exists $A={\rm diag}(A_1,...,A_{\iota_4})\in GL_{2\iota_4}\BF$,
where
$$
A_p=\left(\begin{array}{cc} a_p+b_p&a_p\\
1-a_p-b_p&1-a_p\end{array}\right),\;\;\;
A_q=\left(\begin{array}{cc} 1&0\\ a_q&b_q\end{array}\right) \in
GL_2\BF, \eqno(3.125)$$
for $p\in I_1\cup I_{3,4},\,q\in I_2$,
such that $\a^*_{\{p,\ol p\}}=\a_{\{p,\ol p\}}A_p$ for
$\a\in\G\bs\G_{1,4},\,p\in I_{1,4}$.
\par
Using that $\tau$ is a group isomorphism and applying $\th$ to
$$
[x^{-\si_p},x^\a]=(\a_{\ol p}+\eta_{\ol p}\a_p)x^\a \eqno(3.126)$$
(cf.~(3.3) and (3.91)), by (3.80) and (3.81), we obtain
$$
d_p(\a^*_{\ol p}+\eta_{\ol p}\a^*_p) =\a_{\ol p}+\eta_{\ol p}\a_p
\mbox{ if }\a_{\ol p}+\eta_{\ol
p}\a_p\ne0,\;\a\in\G\bs\G_{1,4},\,p\in I_{1,4} \eqno(3.127)$$
Comparing the coefficients in (3.108), we obtain
$$
(\a_p\g_{\ol p}-\a_{\ol p}\g_p)c_{\si_p+\a+\g} =c_\a
c_\g(\a^*_{p}\g^*_{\ol p}-\a^*_{\ol p}\g^*_{p}) \mbox{ if
}\a_p\g_{\ol p}-\a_{\ol p}\g_p\ne0,\;\a,\g,\a+\g\in\G\bs\G_{1,4}.
\eqno(3.128)$$ Suppose $\a\pm\g\notin\G_{1,4}$. Replacing $\g$ by
$-\g$ in (3.128), and dividing the result from (3.128), we obtain
$$
c_{-\g}c^{-1}_\g=c_{\si_p+\a-\g}c^{-1}_{\si_p+\a+\g}=c_{\a-\g}c^{-1}_{\a+\g}.
\eqno(3.129)$$ In particular, by taking $\g=\si_p+\l_p$ (recall
(3.97)) and replacing $\a$ by $\a+\si_p+\l_p$, we obtain that
$$
c_{-\si_p-\l_p}c_{\a+2\si_p+2\l_p}=c_\a c_{\si_p+\l_p},
\eqno(3.130)$$ holds under some conditions on $\a$ (these
conditions are linear inequalities on $\a_p,\a_{\ol p}$). Setting
$\g=\si_p+2\l_p$ in (3.128) and using (3.130), we obtain
$$
(\a_p(-\eta_{\ol p}+2e_p)-2\a_{\ol p})c^{-1}_{-\si_p-\l_p}
=c^{-1}_{\si_p+\l_p} c_{\si_p+2\l_p}(\a^*_p(-\eta_{\ol
p}+2\l^*_{p,\ol p})-\a^*_{\ol p} (1+2\l^*_{p,p})), \eqno(3.131)$$
holds under some conditions on $\a$. Recall from (3.91) that
$\eta_{\ol p}=0$ if $p\in I_2$ and $\eta_{\ol p}=-1$ otherwise.
Noting that when $p$ is fixed, all coefficients (such as $\l^*_{p,p}$)
of $\a_p,\a_{\ol
p},\a^*_p,\a^*_{\ol p}$ appearing in (3.127) and (3.131) are
constant. {}From (3.127) and (3.131), using (3.99), we can solve
$\a^*_p,\a^*_{\ol p}$ as linear combinations of $\a_p,\a_{\ol p}$
with the coefficient matrices as required in the claim (i.e., as
shown in (3.125)); furthermore, we have $b_p=d_p^{-1}$.
Since $\tau$ is a group isomorphism, the condition
on $\a$ can be removed, i.e., the claim holds for all $\a\in\G$.
\par
{\bf Claim 5}. In (3.125), $a_p=0$ if $p\in I_3$.
\par
Let $p\in I_3$. We write
$\th(t_p)=bt'_p+\sum_{(0,\es_p)\ne(\b,\v j)\in\G'\times\JJ'}b_{\b,\v j}
x^{\b,\v j}$ for some $b,b_{\b,\v j}\in\mbb{F}$. Then we have
$$
\a_{\ol
p}c_{\a+\si_p}x'{}^{\a^*+\si_p}=\th([t_p,x^\a])=bc_\a\a^*_{\ol p}
x'{}^{\a^*+\si_p}+... \eqno(3.132)$$ for $\a\in\G\bs\G_{1,4}$,
where the missed terms do not contain $x'{}^{\a^*+\si_p}$. Thus by
(3.125), we have
$$
\a_{\ol p}c_{\a+\si_p}=bc_\a(a_p\a_p+(1-a_p)\a_{\ol p}).
\eqno(3.133)$$ Hence $b\ne0$. Take $0\ne\a\in\mbb{F}\es_p\cap\G$
(then $\a\notin\G_{1,4}$), we obtain $a_p=0$. This also proves
(3.43) and (3.44).
\par
{\bf Claim 6}. Denote $\si=\sum_{p\in I_{1,4}}\si_p$. For any
$\a\in\G$ with $\a_{_{\sc J_{1,4}}}\ne\si$, we have
$\th(x^\a)=c_\a x'{}^{\a^*}$ for some $c_\a\in\mbb{F}^\times.$
\par
Assume that $\a\in\G_{1,4}$ with $\a_{_{\sc J_{1,4}}}\ne\si$. Then
by (2.26), we can always choose $\b=\b_p\es_p+\b_{\ol p}\es_{\ol
p}\in\G\bs\G_{1,4}$ for some $p\in I_{1,4}$, such that
$$
a=\b_p(\a_{\ol p}-\b_{\ol p}+\eta_p) -\b_{\ol p}(\a_p-\b_p-1)\ne0.
\eqno(3.134)$$
 Then $\b,\a-\b-\si_p\notin\G_{1,4}$ and
$\b^*\in\wt\G'_p$ (where $\wt\G'_p$ is a similar notation as in
(3.107)). We have
$$
\th(x^\a)=a^{-1}\th([x^\b,x^{\a-\b-\si_p}]) =a^{-1}c_\b
c_{\a-\b-\si_p}[x'{}^{\b^*},x'{}^{\a^*-\b^*-\si_p}]
\in\mbb{F}x'{}^{\a^*}. \eqno(3.135)$$
\par
By (3.7) and Lemma 3.1, we have
$$
\th(t_p)\in\HH'{}^F \subset{\rm span}(H'_1\cup H'_3) =\sum_{q\in
I'_{1,4}}\mbb{F}x'{}^{-\si_q}+\sum_{r\in\ol I'_{5,6}}
\mbb{F}t'_r+H'_3, \eqno(3.136)$$ for $p\in\ol I_{5,6}$. Thus,
using notations (2.22) and (3.76), we have (also recall notations
(3.60) and (3.61))
$$
\th^{-1}(\ol t'_{\ol I'_{5,6}})\equiv
(x{}^{-\si})_{I_{1,4}}F_1+\ol t_{\ol I_{5,6}}F_2 \,\ ({\sc\,}{\rm
mod\,\,}H_3{\sc\,}), \eqno(3.137)$$ for some
$$
F_1=(a_{p,q})_{p\in I_{1,4},\,q\in \ol I'_{5,6}} \in
M_{\iota_4\times(\ell'_5+\ell'_6)}, \eqno(3.138)$$
$$
F_2=(b_{p,q})_{p\in \ol I_{5,6},\,q\in \ol I'_{5,6}}\in
GL_{\ell_5+\ell_6}, \eqno(3.139)$$ (in particular
$\ell_5+\ell_6=\ell'_5+\ell'_6$).
\par
{\bf Claim 7}. We have
$$
a_{p,q}=0\mbox{ \  if \ }p\in I_{3,4},\,q\in\ol I'_5,
\eqno(3.140)$$
$$
b_{p,q}=0\mbox{ \ if \ }p\notin I_5,\,q\in\ol I'_5, \eqno(3.141)$$
which implies $(\ell_5,\ell_6)=(\ell'_5,\ell'_6)$ and $I_i=I'_i$
for $i=5,6$.
\par
Note that the center of $\MM$ is $\CC=\{x^\a\in\G\,|\,\a=
\a_{_{\sc I_{5,6}}}\}$. Denote the centralizer
$C_\HH(\CC)=\{u\in\HH\,|\,[u,\CC]=0\}$. It is straightforward to
check that
$$
\{t_p\,|\,p\in I_{2,3}\cup J_4\cup I_6\}\subset C_\HH(\CC)\subset
{\rm span}\{x^{\a,\v i}\,|\,\v i_{_{\sc \ol I_{5,6}}}=0\}.
\eqno(3.142)$$
For $p\in I_{1,4}$, (3.142) implies that
$\ad_{x^{-\si_p}}|_{C_\HH(\CC)}$ is semi-simple if and only if
$p\in I_{1,2}$, and ${\rm ad}_{t_q}|_{C_\HH(C)}$ is semi-simple
for $q\in\ol I_5$ by (3.142) and is not semi-simple for $q\in\ol
I_6$. Moreover, by (3.1), for $p\in\ol I_{5,6}$, ${\rm ad}_{t_p}$
is semi-simple if and only if $p\in\ol I_5$. We obtain the claim.
\par
By (3.140) and (3.141), we can write $F_1$ and $F_2$ is the forms
$$
F_1=\pmatrix{B_{1,5}&B_{1,6}\cr
B_{2,5}&B_{2,6}\cr0&B_{3,6}\cr},\;\;\;
F_2=\pmatrix{B_{5,5}&B_{5,6}\cr0&B_{6,6}\cr}, \eqno(3.143)$$ such
that all $B_{i,j}$ have the forms in (3.46) and (3.48).
\par
For any $\a\in\G$, we denote
$$
\wh\a=(\a_{\ol1}+\eta_{\ol1}\a_1,...,\a_{\ol\iota_4}+\eta_{\ol\iota_4}
\a_{\iota_4}) \in\mbb{F}^{\iota_4}, \eqno(3.144)$$ (cf.~(3.91)).
For $\a\in\G\bs\G_{1,4}$, applying $\th^{-1}$ to
$$
\a^*_{I_{5,6}} c_\a^{-1}x'{}^{\a^*}= [\ol
t'_{I_{5,6}},c_\a^{-1}x'^{\a^*}], \eqno(3.145))$$ (cf.~(3.4)),
using (3.137), and noting that $[H_3,\MM]=0$, we obtain
$$
\a^*_{I_{5,6}}x^{\a}= [(x{}^{-\si})_{I_{1,4}}F_1+\ol t_{\ol
I_{5,6}}F_2,x^\a] =(\wh\a F_1+\a_{_{\sc I_{5,6}}}F_2)x^\a,
\eqno(3.146)$$ that is
$$
\a^*_{I_{5,6}}=\wh\a F_1+\a_{_{\sc I_{5,6}}}F_2, \eqno(3.147)$$
holds for all $\a\in\G\bs\G_{1,4}$ and so holds for all $\a\in\G$
since $\tau:\a\mapsto\a^*$ is an isomorphism. {}From this and
(3.143), we obtain formulas (3.43)-(3.48) (cf.~(3.49)).
\par
{\bf Claim 8}. $\ell_7=\ell'_7$.
\par
Observe from (3.9) that
$$
H_3=C_\HH(\MM)\ \mbox{ ( the centralizer of \ $\MM$)},
\eqno(3.148)$$
$$
{\rm span}\{x^{\a,\v i}\in H_3\,|\,\v i_{_{\sc J_7}}=0\} =C(H_3)\
\mbox{ ( the center of \ $H_3$)}. \eqno(3.149)$$ By exchanging
$\HH$ with $\HH'$ if necessary, we can suppose $\ell_7\le\ell'_7$.
As in the proof of sufficiency, we can construct an embedding
$\ol\th:\HH\rar\HH'$ such that
$$
\ol\th(x^\a)=\th(x^\a),\;\;\; \ol\th(\ol t_{\ol
I_{5,6}})\equiv\th(\ol t_{\ol I_{5,6}}) \,\ ({\sc\,}{\rm
mod\,\,}H'_3{\sc\,}), \eqno(3.150)$$ (cf.~Claim 6 and (3.137),
note that using (3.137), we can now obtain that Claim 6 holds for
all $\a\in\G$ if $\ell_5+\ell_6\ne0$). Thus by identifying $\HH$
with $\ol\th(\HH)$, we can assume that $\HH$ is a subalgebra of
$\HH'$ such that there exists an isomorphism $\th$ satisfying
$$
\th(x^\a)=x^\a,\;\;\th(t_{\ol p})\equiv t_{\ol p} \,\ ({\sc\,}{\rm
mod\,\,}H_3'{\sc\,}) \;\;\;\for\;\;\;\a\in\G,p\in I_{5,6}.
\eqno(3.151)$$ By restricting $\th$ to $H_3$, we want to prove
$$
\th(t_p)=t_p+c_p \;\;\;\for\;\;\;p\in I_6\mbox{ \,and some
\,}c_p\in\mbb{F}, \eqno(3.152)$$
$$
\th(x^{\a,\v i} t^{\v j})= x^\a\prod_{p\in
I_6}(\th(t_p))^{i_p}\prod_{q\in J_7}(\th(t_q))^{j_q}
\;\;\;\for\;\;\;\a=\a_{I_{5,6}},\v i=\v i_{I_6},\v j=\v j_{J_7}.
\eqno(3.153)$$ To prove (3.152), first by (3.149), we have
$c_p=\th(t_p)-t_p\in C(H'_3)$. Then by (3.151), we have
$$
[t_{\ol q},c_p]= \th([\th^{-1}(t_{\ol q}),t_p])-[t_{\ol q},t_p]
=0, \eqno(3.154)$$ where the second equality follows from the fact
that $\th^{-1}(t_{\ol q})=t_{\ol q}\,(\,{\rm mod\,}H_3)$ and
$[H_3,t_p]=0$. {}From (3.154), we obtain that $c_p\in\mbb{F}$. Thus
we have (3.152). Similarly, we have
$$
\th(x^{\a,\es_{p}})=x^\a(t_p+c_{\a,p}) \;\;\;\for\;\;\;p\in
I_6\mbox{ and some }c_{\a,p}\in\mbb{F}. \eqno(3.155)$$ By
considering $\th([x^\a,t_pt_{\ol p}])=[\th(x^\a),\th(t_pt_{\ol p})]$,
we see that $c_{\a,p}=c_p$, and we obtain
$$
\th(t_pt_{\ol p})=(t_p+c_p)t_{\ol p}+u_p \;\;\;\for\;\;\;p\in I_6
\mbox{ and some }u_p\in C_{\HH'}(C(H'_3)). \eqno(3.156)$$ {}From
this and (3.152), we can deduce
$$
\th(x^{\a,\es_{p}})=x^\a(t_p+c_p)\;\;\;\for\;\;\;p\in I_6.
\eqno(3.157)$$ Similar to (3.156), we have
$$
\th(x^{-\si_p,\es_{p}+\es_{\ol p}})= x^{-\si_p,\es_{\ol
p}}(t_p+c_p)+u'_p \;\;\;\for\;\;\;p\in I_6 \mbox{ and some
}u'_p\in C_{\HH'}(C(H'_3)). \eqno(3.158)$$ Now from  (3.152),
(3.155)-(3.158), we can obtain (3.153) by induction on $|\v i|$ in
case $\v j=0$.
\par
Assume that (3.153) holds for all $\v j$ with $|\v j|<n$, where
$n\ge1$. We denote by $A_{\a,\v i,\v j}$ the difference between
the left-hand side and the right-hand side of (3.153). Then the
inductive assumption says that $A_{\a,\v i,\v j}=0$ if $|\v j|<n$.
Now suppose $|\v j|=n$. Say $j_r\ge1$ for some $r\in I_7$ (the
proof is similar if $r\in\ol I_7$). Let $\v k=\v
j-\es_{r}+\es_{\ol r}$. Then we have
$$
\begin{array}{ll}
[\th(t_r),A_{\a,\v i,\v k}]
\!\!\!\!&\dis=\th([t_r,x^{\a,\v i}t^{\v k}])
-\th([t_r,\th^{-1}(x^\a)])
\prod_{p\in I_6}(\th(t_p))^{i_p}\prod_{q\in J_7}(\th(t_q))^{k_q}
\vs{4pt}\\
&\dis\ \ \ \ -x^\a[\th(t_r),\prod_{p\in I_6}(\th(t_p))^{i_p}\prod_{q\in J_7}
(\th(t_q))^{k_q}]
\vs{4pt}\\
&\dis=(j_{\ol r}+1)(
\th(x^{\a,\v i}t^{\v j-\es_{r}})-
x^\a\prod_{p\in I_6}(\th(t_p))^{i_p}\prod_{q\in J_7}
(\th(t_q))^{j_q-\d_{q,r}})
\vs{4pt}\\
&\dis=(j_{\ol r}+1)A_{\a,\v i,\v j-\es_{r}}=0,
\end{array}
\eqno(3.159)$$ where the first equality follows from (1.1), the
second equality follows from (1.1) and (3.151). By (1.1) and
(3.159), we obtain
$$
[\th(t^2_r),A_{\a,\v i,\v k}]= \th([t^2_r,\th^{-1}(A_{\a,\v i,\v
k})]) =2\th(t_r[t_r,\th^{-1}(A_{\a,\v i,\v k})]) =0.\eqno(3.160)$$
On the other hand, exactly similar to (3.159), we have
$$
[\th(t^2_r),A_{\a,\v i,\v k}]= 2(j_{\ol r}+1)A_{\a,\v i,\v j}.
\eqno(3.161)$$ Now (3.160) and (3.161) show that $A_{\a,\v i,\v
j}=0$. This proves (3.153). By (3.152), (3.153) and by identifying
$C(H_3)$ with $C(H'_3)$ using the isomorphism, we see that $\th$
is an associative algebra isomorphism $H_3\to H'_3$ over the
domain ring $C(H_3)$. {}From this we obtain $\ell_7=\ell'_7$ since
$2\ell_7$ is the transcendental degree of $H_3$ over the domain
ring $C(H_3)$. This completes the proof of Theorem 3.6. \qed\par
\vs{5pt}\par\ni {\bf4. \ Derivations} \vs{-1pt}\par\ni In this
section, we shall determine the structure of the derivation
algebra of the Hamiltonian Lie algebra $\HH=\HH(\v\ell,\G)$. As
pointed in [F], the significance of derivations for Lie theory
primarily resides in their affinity to low dimensional cohomology
groups, their determination therefore frequently affords insight
into structural features of Lie algebras which do not figure
prominently in the defining properties. Some general results
concerning derivations of graded Lie algebras were established in
[F]. However in our case the algebras are in general nongraded,
the results in [F] can not be applied to our case here. Thus we
try a different method to determine derivations of the Hamiltonian
Lie algebras $\HH$. Our method is also different from that used in
[OZ].
\par
Recall that a {\it derivation} $d$ of the Lie algebra $\HH$ is a linear
transformation on $\HH$ such that
$$
d([u_1,u_2])=[d(u_1),u_2]+[u_1,d(u_2)]\ \ \for\ \ u_1,u_2\in\HH.
\eqno(4.1)$$
Denote by $\der\HH$ the space of the derivations of $\HH$,
which is a Lie algebra. Moreover, $\ad_\HH$ is an ideal.
Elements in $\ad_\HH$ are called
{\it inner derivations}, while elements in $\der\HH\bs\ad_\HH$ are
called {\it outer derivations}.
\par
We can embed $\HH$ into a larger Lie algebra $\wt\HH$ such that
$\wt\HH$ has a basis $\{x^{\a,\v i}\,|\,(\a,\v i)\in
\G\times\mbb{N}^{2\iota_7}\}$ (i.e., in $\wt\HH$, we replace $\JJ$
by $\mbb{N}^{2\iota_7}$, cf.~(2.27), and we have (3.1) with the
last three summands running over $p\in I_{1,6},\,p\in I_{1,4}$ and
$p\in I$ respectively). Then for $p\in J_1\cup\ol I_{2,3}\cup
I_5$, clearly, $t_p\notin\HH$, but $[t_p,\HH]\subset\HH$. Thus
$$
d_p=\ad_{t_p}|_\HH\;\;\;\for\;\;\;p\in J_1\cup\ol I_{2,3}\cup I_5,
\eqno(4.2)$$
defines an outer derivation of $\HH$. For
$p\in I_{2,3}\cup J_4\cup\ol I_5\cup J_{6,7}$, obviously, $\ptl_{t_p}$ is a
derivation of $\HH$ (cf.~(2.27), (2.35) and (3.1)). For $p\in J$,
we define $\SGN(p)=1$ if $p\in I$ and $\SGN(p)=-1$ if $p\in\ol I$.
Then
$$
\ptl_{t_p}=\SGN(p)\ad_{t_{\ol p}}\;\;\;\for\;\;\;p\in\ol I_5\cup J_{6,7}.
\eqno(4.3)$$
Define $d_0(x^{\a,\v i})=(\sum_{p\in I_{1,4}}\a_p+1)x^{\a,\v i}$ for
$(\a,\v i)\in\G\times\JJ$. It is straightforward to verify that $d_0$ is an
outer derivation of $\HH$.
Denote $\si=\sum_{p\in I_{1,4}}\si_p$.
If $\iota_7=\ell_1$, then $\HH=[\HH,\HH]+\mbb{F}x^\si$, and we can define
an outer derivation $d'_0$ by setting
$$
d'_0([\HH,\HH])=0,\;\;\;d'_0(x^\si)=1_{\HH}.
\eqno(4.4)$$
If $\iota_7\ne\ell_1$, we set $d'_0=0$.
\par
We denote by $\HOM$ the set of group homomorphisms
$\mu:\G\to\mbb{F}$ such that $\mu(\si_p)=0$ for $p\in I_{1,4}$.
For $\mu\in\HOM$, we define a linear transformation $d_\mu$ on
$\HH$ by
$$
d_{\mu}(x^{\a,\v i})=\mu(\a)x^{\a,\v i}\;\;\;\for\;\;\;(\a,\v i)
\in\G\times\JJ.
\eqno(4.5)$$
Clearly, by (3.1), $d_\mu$ is a derivation of $\HH$. We identify
$\HOM$ with a subspace of $\der\HH$ by $\mu\mapsto d_\mu$. For
$p\in I_{1,6}$, we define $\mu_p\in\HOM$ by
$$
\mu_p(\a)=
\left\{\matrix{\a_{\ol p}+\eta_{\ol p}\a_p\hfill&\mbox{if \ \ }p\in I_{1,4},
\vs{4pt}\hfill\cr
\a_p\hfill&\mbox{if \ \ }p\in I_{5,6},\hfill\cr}\right.
\eqno(4.6)$$
for $\a\in\G$ (cf.~(3.91)). By (3.3) and (3.4), we have
$$
\ad_{x^{-\si_p}}=\left\{\matrix{
d_{\mu_p}\hfill&\mbox{if \ \ \ }p\in I_{1,2},
\vs{4pt}\hfill\cr
d_{\mu_p}+\ptl_{t_p}\hfill&\mbox{if \ \ \ }p\in I_3,
\vs{4pt}\hfill\cr
d_{\mu_p}+\ptl_{t_p}-\ptl_{t_{\ol p}}
\hfill&\mbox{if \ \ \ }p\in I_4,
\hfill\cr}\right.
\eqno(4.7)$$
$$
\ad_{t_{\ol q}}=\left\{\matrix{ -d_{\mu_q}\hfill&\mbox{if \ \ \
}q\in I_5, \vs{4pt}\hfill\cr -d_{\mu_q}-\ptl_{t_q}\hfill&\mbox{if
\ \ \ }q\in I_6.
 \hfill\cr}\right. \eqno(4.8)$$ We fix a subspace
$\hom$ of $\HOM$ such that
$$
\HOM=\hom\oplus{\rm span}\{\mu_p\,|\,p\in I_{1,6}\},
\eqno(4.9)$$
is a direct sum as vector spaces.
Since $\ad_{\mbb{F}}=0$, we set
$\HH^*={\rm span}\{x^{\a,\v i}\,|\,(0,0)\ne(\a,\v i)\in\G\times\JJ\}$.
\par\ni
{\bf Theorem 4.1}.
{\it The derivation algebra $\der\HH$ is spanned by
$$
d'_0,\;d_p,\;\ptl_{t_q},\;d_\mu,\;\ad_{\HH^*}\;\;\;\for\;\;\;p\in\{0\}\cup
J_1\cup\ol I_{2,3}\cup I_5,\,q\in I_{2,3}\cup J_4,\,\mu\in\hom.
\eqno(4.10)$$
Furthermore, we have the following vector space decomposition as a direct
sum of subspaces:
$$
\der\HH=((\mbb{F}d'_0+\sum_{p\in\{0\}\cup J_1\cup\ol I_{2,3}\cup I_5}
\mbb{F}d_p)
\oplus\sum_{q\in I_{2,3}\cup J_4}\mbb{F}\ptl_{t_q}
\oplus\hom\,)\oplus\ad_{\HH^*}.
\eqno(4.11)$$
In particular, all derivations of the classical Hamiltonian Lie algebras
$\HH(\ell)$ (cf.~(1.3)) are inner.
}
\par\ni
{\it Proof.}
First note that in [OZ], $d_0$ was written as a derivation of the form $d_\mu$
with $\mu$ satisfying $\mu(\si_p)=\mu(\si_1)$ for $p\in I_{1,4}$.
Let $d\in\der\HH$ and let $D$ be the subspace of $\der\HH$ spanned by
the elements in (4.10).
Note that $D\supset\HOM$ by (4.7) and (4.8).
We shall prove that after a number of steps in each of
which $d$ is replaced by $d-d'$ for some $d'\in D$ the $0$ derivation is
obtained and thus proving that $d\in D$. This will be done by a number of
claims.
\par
{\bf Claim 1}.
We can suppose (i) $d(1)=0$, (ii) $d(x^{-\si_p})=0$ for $p\in I_{3,4}$,
(iii) $d(x^{\es_q})=d(t_r)=0$ for $q\in I_{5,6},\,r\in\ol I_6\cup J_7$.
\par
By replacing $d$ by $d-d(1)d_0$, we can suppose $d(1)=0$.
For any $(\a,\v i)\in\G\times\JJ$, we write
$$
d(x^{\a,\v i})=\sum_{(\b,\v j)\in M_{\a,\v i}}
c_{\a,\v i}^{(\b,\v j)}x^{\a+\b,\v j}\;\;\;\mbox{for some}\;\;\;c_{\a,\v i}^{(\b,\v j)}\in\mbb{F},
\mbox{ where}
\eqno(4.12)$$
$$
M_{\a,\v i}=\{(\b,\v j)\in\G\times\JJ\,|\,c_{\a,\v i}^{(\b,\v j)}\ne0\},
\eqno(4.13)$$
is a finite set. We set $c_{\a,\v i}^{(\b,\v j)}=0$
if $(\b,\v j)\notin M_{\a,\v i}$.
We shall denote $M_{\a,0}$ simply by $M_\a$.
Using inductive assumption, suppose we have proved that
$d(x^{-\si_r})=0$ for $r\in I_{3,4}$ and $r<p$.
Let $(\b,\v j)\in M_{-\si_p}$. Using (3.3), one can deduce by induction on
$|\v j|$ that
$$
x^{-\si_p+\b,\v j}=[u_{\b,\v j},x^{-\si_p}]\;\;\; \mbox{for
some}\;\;\;u_{\b,\v j}\in\HH, \eqno(4.14)$$ such that $u_{\b,\v
j}$ has the following form
$$
u_{\b,\v j}=\sum_{k,\ell\in\mb{Z}}
b_{k,\ell}x^{-\si_p+\b,\v j+k\es_p+\ell\es_{\ol p}}\;\;\;\mbox{for some}\;\;\;
b_{k,\ell}\in\mbb{F},
\eqno(4.15)$$
(recall convention (2.34)).
Thus we can take
$$
u=\sum_{(\b,\v j)\in M_{-\si_p}}c_{-\si_p,0}^{(\b,\v j)}
u_{\b,\v j}\in\HH
\mbox{ \ such that \ }(d-\ad_u)(x^{-\si_p})=0,
\eqno(4.16)$$
Applying $d$ to $[x^{-\si_r}, x^{-\si_p}]=0$, we obtain
$$
\sum_{(\b,\v j)\in M_{-\si_p}}c_{-\si_p,0}^{(\b,\v j)}
[x^{-\si_r},x^{-\si_p+\b,\v j}]=0\;\;\;\for\;\;\;r\in I_{3,4},\,r<p,
\eqno(4.17)$$
i.e.,
$$
c_{-\si_p,0}^{(\b,\v j)}(\b_r-\b_{\ol r})
-c_{-\si_p,0}^{(\b,\v j+\es_{\ol r})}(j_{\ol r}+1)
+c_{-\si_p,0}^{(\b,\v j+\es_r)}(j_r+1)=0,
\eqno(4.18)$$
for $r\in I_{3,4},\,r<p,$ from this and by induction on $j_r+j_{\ol r}$
ranging from  ${\rm max}\{k_r+k_{\ol r}\,|\,(\b,\v k)\in M_{-\si_p}\}$
down to zero, we obtain
$$
\b_r=\b_{\ol r}, \,j_r=j_{\ol r}=0\;\;\;\for\;\;\;(\b,\v j)\in
M_{-\si_p},\, r\in I_{3,4},\,r<p. \eqno(4.19)$$ Then (4.15),
(4.16) and (4.19) show that $\ad_u(x^{-\si_r})=0$ for $r\in
I_{3,4},\,r<p$. Thus if we replace $d$ by $d-\ad_u$, we have
$d(x^{-\si_r})=0$ for $r\in I_{3,4},\,r\le p$. This proves Claim
1(ii). Note that for $v=x^{\es_q},\,q\in I_{5,6},$ or
$v=t_r,\,r\in\ol I_6\cup J_7$, we have $\ad_v(\HH)=\HH$. Thus
similar to the above proof, we have Claim 1(iii).
\par
Note that for $(\b,\v j)\in\G\times\JJ$, by (3.3) and (3.4), we have
$$
(\b_{\ol p}+\eta_{\ol p}\b_p)x^{-\si_p+\b,\v j}
=[x^{-\si_p+\b,\v j},x^{-\si_p}]\;\;\;\for\;\;\;p\in I_{1,2},
\eqno(4.20)$$
$$
(-1+\b_p)e_p x^{\l_p+\b,\v j}=[x^{-\si_p+\b,\v j},x^{\l_p}]
+j_px^{\l_p+\b,\v j-\es_p}\;\;\;\for\;\;\;p\in I_{1,4},
\eqno(4.21)$$ (recall notations $\l_p,\,p\in I_{1,4}$ in (3.97)),
and
$$
\b_px^{\b,\v j}=[x^{\b,\v j},t_{\ol p}]\;\;\;\for\;\;\;p\in I_5.
\eqno(4.22)$$
\par
{\bf Claim 2}.
By replacing $d$ by $d-d'$ for some $d'\in D$, we can suppose
$$
\b_{\ol p}+\eta_{\ol p}\b_p=0\;\;\;\for\;\;\;(\b,\v j)\in M_{-\si_p},
\,p\in I_{1,2},
\eqno(4.23)$$
$$
\b_p=1\;\;\;\for\;\;\;(\b,\v j)\in M_{\l_p},\,
p\in I_{1,4},
\eqno(4.24)$$
$$
\b_p=0\;\;\;\for\;\;\;(\b,\v j)\in M_{0,\es_{\ol p}},\,
p\in I_5.
\eqno(4.25)$$
\par
The proof of (4.23) is similar to that of Claim 1. To prove (4.24), suppose
we have proved
$$
\b_r=1\;\;\;\for\;\;\;(\b,\v j)\in M_{\l_r},\,
i\in I_{1,4},\;r<p.
\eqno(4.26)$$
To see how the proof works, for
simplicity, we assume that $p\in I_1$ (the proof for $p\in I_{2,4}$
is exactly similar). Then the second term on the right-hand
side of (4.21) vanishes.
Let
$$
u=\sum_{(\b,\v j)\in M_{\l_p},\,\b_p\ne1}c_{\l_p,0}^{(\b,\v j)}
((-1+\b_p)e_p)^{-1}x^{-\si_p+\b}.
\eqno(4.27)$$
Then by replacing $d$ by $d-\ad_u$, from (4.21), we see that (4.24) holds for
$p$. We want to prove that after this replacement, Claim 1, (4.23) and
(4.26) still hold. It suffices to prove
$$
[u,x^{-\si_q}]=[u,x^{\es_{q'}}]=[u,t_{q''}]=[u,x^{\l_r}]=0,
\eqno(4.28)$$
for $q\in I_{1,4},\,q'\in I_{5,6},\,q''\in\ol I_6\cup J_7,\,r\in I_{1,4},\,
r<p.$
\par
We have
$$
\matrix{\dis
-e_p\sum_{(\b,\v j)\in M_{\l_p}}c_{\l_p,0}^{(\b,\v j)}x^{\l_p+\b,\v j}
=-e_pd(x^{\l_p})
=d([x^{-\si_p},x^{\l_p}]
\vs{4pt}\hfill\cr\dis\ \ \
=
\sum_{(\b,\v j)\in M_{-\si_p}}c_{-\si_p,0}^{(\b,\v j)}(-1+\b_p)
x^{\l_p+\b,\v j}
+
\sum_{(\b,\v j)\in M_{\l_p}}c_{\l_p,0}^{(\b,\v j)}(\b_p-\b_{\ol p}-e_p)
x^{\l_p+\b,\v j}.
\hfill\cr}
\eqno(4.29)$$
This gives
$$
\b_p-\b_{\ol p}=
(c_{\l_p,0}^{(\b,\v j)})^{-1}
(\b_p-1)c_{-\si_p,0}^{(\b,\v j)}
\;\;\;\for\;\;\;(\b,\v j)\in M_{\l_p}.
\eqno(4.30)$$
If $(\b,\v j)\notin M_{-\si_p}$, then the right-hand side of (4.30) is zero;
on the other hand, if $(\b,\v j)\in M_{-\si_p}$, then (4.23) gives
$\b_p-\b_{\ol p}=0$. In any case, we have $\b_p-\b_{\ol p}=0$
for $(\b,\v j)\in M_{\l_p}$. Thus by (4.27),
$$
[u,x^{-\si_p}]=
\sum_{(\b,\v j)\in M_{\l_p},\b_p\ne1}c_{\l_p,0}^{(\b,\v j)}
((-1+\b_p)e_p)^{-1}
(\b_p-\b_{\ol p})x^{-\si_p+\b}=0.
\eqno(4.31)$$
Similarly, we can prove other equations in (4.28). This proves (4.24).
Similarly, we have (4.25).
\par
{\bf Claim 3}. By replacing $d$ by
$d-\sum_{p\in J_1\cup\ol I_{2,3}\cup I_5}a_pd_p-d_\mu$
for some $a_p\in\mbb{F}$ and some $\mu\in\HOM$,
we can suppose $d(x^{-\si_p})=d(x^{\l_q})=d(t_r)=0$ for
$p\in I_{1,2},\,q\in I_{1,4},\,r\in I_5.$
\par
Again for simplicity, we prove that after some replacement,
$d(x^{-\si_p})=d(x^{\l_p})=0$ for $p\in I_1$. Defining
$\mu\in\HOM$ by $\mu(\a)=c_{\l_p,0}^{(0,0)}e_p^{-1}(\a_{\ol
p}-\a_p)$, and by replacing $d$ by $d-d_\mu$, we obtain
$c_{\l_p,0}^{(0,0)}=0$ (recall (4.12) that $c_{\a,\v i}^{(\b,\v
j)}$ is the coefficient of $x^{\a+\b,\v j}$, not that of $x^{\b,\v
j}$). Obviously, this replacement does not affect the result we
have obtained so far. Recalling the definition of $d_p$ in (4.2),
we have
$$
d_p(x^{-\si_p})=[t_p,x^{-\si_p}]=-1,\,\;\;d_{\ol p}(x^{-\si_p})=1,\,\;\;
d_p(x^{\l_p})=e_px^{\si_p+\l_p},\,\;\;d_{\ol p}(x^{\l_p})=0.
\eqno(4.32)$$
Thus by replacing $d$ by $d-a_pd_p-a_{\ol p}d_{\ol p}$ for some $a_p,a_{\ol p}
\in\mbb{F}$, we can suppose
$$
c_{-\si_p,0}^{(\si_p,0)}=c_{\l_p,0}^{(\si_p,\l_p)}=c_{\l_p,0}^{(0,0)}=0,
\eqno(4.33)$$ Note again that the replacement does not affect the
results we have obtained so far.
\par
Let $q\in I_1,\,q\ne p$. We have
$$
0=d([x^{-\si_p},x^{-\si_q}])=
\sum(c_{-\si_p,0}^{(\b,\v j)}(\b_{\ol q}-\b_q)
+c_{-\si_q,0}^{(\b+\si_q-\si_p,\v j)}(\b_p-\b_{\ol p}))x^{-\si_p+\b,\v j},
\eqno(4.34)$$
$$
0=d([x^{-\si_p},x^{\l_q}])=
\sum(c_{-\si_p,0}^{(\b,\v j)}\b_qe_q+c_{\l_q,0}^{(\b+\si_q-\si_p,\v j)}
(\b_p-\b_{\ol p}))x^{-\si_p+\si_q+\b+\l_q,\v j}.
\eqno(4.35)$$
Now (4.23), (4.34) and (4.35) show that $\b_{\ol q}=\b_q=0$
if $(\b,\v j)\in M_{-\si_p}$. Similarly,
we can prove $\b_r=j_r=0$ for all $r\in J,\,r\ne p,\ol p$ if
$(\b,\v j)\in M_{-\si_p}$. This and (4.30) show that $(\b,\v j)=(\si_p,0)$
if $(\b,\v j)\in M_{-\si_p}$. But $(\si_p,0)\notin M_{-\si_p}$ by (4.33),
i.e., $M_{-\si_p}=\emptyset$. Thus $d(x^{-\si_p})=0$. Similarly
$d(x^{\l_p})=0$. Analogously, we can obtain other results of Claim 3.
\par
{\bf Claim 4}. We can suppose $d=0$.
\par
Note that $x^\a$ is a common eigenvector for the elements of the set
$$
A=\{x^{-\si_p},\,x^{\es_q},\,t_r\,|\,p\in I_{1,4},\,q\in I_{5,6},\,r\in
\ol I_{5,6}\cup J_7\}.
\eqno(4.36)$$
Since $d(A)=0$, $d(x^\a)$ is also a common eigenvector for the elements of
$A$. {}From this and Lemma 3.2, we obtain
$$
\eta_{\ol p}\b_p+\b_{\ol p}=\b_q=0,\,\;\v j=0\;\;\;
\for\;\;\;p\in I_{1,4},\,q\in I_{5,6}\;\mbox{ and }\;(\b,\v j)\in M_\a.
\eqno(4.37)$$
For simplicity, we denote $c_\a^{(\b)}=c_{\a,0}^{(\b,0)}$.
We want to prove
$$
d(x^\a)=m_\a x^\a\mbox{ \ for \ $\a\in\G$ \ and some \ }
m_\a\in\mbb{F}, \eqno(4.38)$$ i.e., $M_\a$ is either empty or a
singleton $\{(0,0)\}$. Thus assume that
$$
\b_p\ne0\mbox{ \ for some \ }(\b,\v j)\in M_\a,\;p\in I_{1,4},\;\a\in\G.
\eqno(4.39)$$
\par
For convenience, we again suppose $p\in I_1$.
Denote $\G_p=(\mbb{F}\es_p+
\mbb{F}\es_{\ol p})\cap\G$ as in (3.102),
and set $\HH_p={\rm span}\{x^\a\,|\,\a\in\G_p\}$.
We have
$$
d(x^\a)\in\HH_p\;\;\;\for\;\;\;\a\in\G_p, \eqno(4.40)$$ by using
the fact that $x^\a$ commutes with elements of $A$ except possibly
$x^{-\si_p},x^{\l_p}$. By (4.37) and by
$$
0=d([x^{-\si_p-\l_p},x^{\l_p}])=
[d(x^{-\si_p-\l_p}),x^{\l_p}]=
\sum c_{-\si_p-\l_p}^{(\b)}e_px^\b,
\eqno(4.41)$$
and by (4.40), we obtain
$$
d(x^{-\si_p-\l_p})=a_px^{-\l_p},\;\;\;
d(x^{\l_p-\si_p})=-a_px^{\l_p}\mbox{ \ for some \ }a_p\in\mbb{F},
\eqno(4.42)$$
where the second equation is obtained from
$d([x^{-\si_p-\l_p},x^{\l_p-\si_p}])=0$.
Applying $d$ to
$$
[x^{-2\si_p},x^{\l_p}]=-2e_px^{-\si_p+\l_p},\;\;\;\;
[x^{-2\si_p},x^{-\si_p-\l_p}]=2e_px^{-2\si_p-\l_p},
\eqno(4.43)$$
$$
[x^{-2\si_p-\l_p},x^{\l_p-\si_p}]=-3e_px^{-2\si_p},
\eqno(4.44)$$
we obtain respectively
$$
d(x^{-2\si_p})=-2a_px^{-\si_p},\;\;
d(x^{-2\si_p-\l_p})=0,\;\;
a_p=0.
\eqno(4.45)$$
Thus all equations in (4.42) and (4.45) are zero. Applying
$d$ to $[x^{\l_p-\si_p},x^{k\l_p}]=-kx^{(k+1)\l_p}$,
using induction on $k$, we obtain
$$
d(x^{k\l_p})=0\;\;\;\;\for\;\;\;k\ge1.
\eqno(4.46)$$
\par
Applying $d$ to
$$
[x^\a,x^{-2\si_p}]=2(\a_{\ol p}-\a_p)x^{\a-\si_p},\;\;
[x^\a,x^{-\si_p-k\l_p}]=(\a_{\ol p}-\a_p-k\a_{\ol p}e_p)x^{\a-k\l_p},
\eqno(4.47)$$
$$
[x^{\a-k\l_p},x^{k\l_p}]=k\a_pe_px^{\a+\si_p},
\eqno(4.48)$$
for $k\ge1$,
using (4.37), we obtain
$$
2(\a_{\ol p}-\a_p)c_{\a}^{(\b)}=2(\a_{\ol p}-\a_p)c_{\a+\si_p}^{(\b)},
\eqno(4.49)$$
$$
(\a_{\ol p}-\a_p-k(\a_p+\b_p)e_p)c_{\a}^{(\b)}=(\a_{\ol p}-\a_p-k\a_pe_p)
c_{\a-k\l_p}^{(\b)},
\eqno(4.50)$$
$$
k(\a_p+\b_p)e_pc_{\a-k\l_p}^{(\b)}=k\a_pe_pc_{\a+\si_p}^{(\b)}.
\eqno(4.51)$$
If $\a_p\ne\a_{\ol p}$, then the above three equations
gives $\b_p=0$,
a contradiction with (4.39). Thus we obtain
$$
\b_p\ne0,(\b,0)\in M_\a\ \Rar\ \a_p=\a_{\ol p}.
\eqno(4.52)$$
Replacing $\a$ by $\a-\si_p$ in (4.51),
it gives
$$
(\a_p-1)c_\a^{(\b)}=(\a_p-1+\b_p)c_{\a-\si_p-\l_p}^{(\b)}.
\eqno(4.53)$$ Noting that for $\a'=\a-\si_p-\l_p$, we have
$\a'_p\ne\a'_{\ol p}$. Assume $(\b,0)\in M_\a$. If $(\b,0)\in
M_{\a'}$, then (4.52) shows that $\b_p=0$, again a contradiction
with (4.39). Thus $(\b,0)\notin M_{\a'}$ and the right-hand side
of (4.53) is zero. This and (4.52) show that $\a_{\ol p}=\a_p=1$.
Note that for $\a''=\a-k\l_p,\,k\ge1$, the relation $\a''_{\ol
p}=\a''_p=1$ does not hold, thus the right-hand side of (4.50) is
zero. We obtain $\a_p+\b_p=0$. Hence
$$
\a_{\ol p}=\a_p=-\b_p=-\b_{\ol p}=1 \;\mbox{ if
}\;\b_p\ne0,\;(\b,0)\in M_\a. \eqno(4.54)$$ If $\a_{_{\sc
J_{1,4}}}\ne\si$ (cf.~(2.22)), say $(\a_q,\a_{\ol q})\ne(1,1)$ for
some $q\in I_1,\,q\ne p$. Suppose $\a_q\ne1$ (the proof is similar
if $\a_{\ol q}\ne1$), then we can write
$$
x^\a = ((\a_q-1)e_q)^{-1}
[x^{\a-\si_q-\l_q+\si_p},x^{\l_q-\si_p}]. \eqno(4.55)$$
 Since for
$\a'=\a-\si_q-\l_q+\si_p$ or $\l_q-\si_p$, the relation
$\a'_p=\a'_{\ol p}=1$ does not holds, we have $\b_p=0$ if
$(\b,0)\in M_{\a'}$. Then applying $d$ to (4.55) gives that
$\b_p=0$ if $(\b,0)\in M_\a$, which contradicts (4.39) again.
Hence
$$
\a_{_{\sc J_{1,4}}}=\si,\mbox{ \ and \ }\b=-\si, \eqno(4.56)$$ by
(4.37) and (4.54). If $\ell_5+\ell_6+\ell_7\ne0$, we can write
$$
x^\a=\left\{\matrix{
(\a_q+k)^{-1}[x^{\a+\si_p+k\es_q},x^{-\si_p-k\es_q,\es_{\ol q}}]
\hfill&\mbox{for \ }q\in I_{5,6},\,k\in\mb{Z},\,\a_q+k\ne0,
\vs{4pt}\hfill\cr
[x^{\a+\si_p,\es_r},x^{-\si_p,\es_{\ol r}}]
\hfill&\mbox{for \ }r\in I_7.
\hfill\cr}\right.
\eqno(4.57)$$
Note that for
$$
(\a',\v i')=(\a+\si_p+k\es_q,0),\,
(-\si_p-k\es_q,\es_{\ol q}),\,(\a+\si_p,\es_r)
\mbox{ \,or \,}(-\si_p,\es_{\ol r}),
\eqno(4.58)$$
the relation $\a'_p=\a'_{\ol p}=1$ does not holds,
one can prove as above that
$\b_p=0$ if $(\b,0)\in M_{\a',\v i'}$.
Then applying $d$ to (4.57) gives that
$\b_p=0$ if $(\b,0)\in M_\a$, which again contradicts (4.39). Hence
$\ell_5+\ell_6+\ell_7=0$. Similarly, one can prove $\ell_2+\ell_3+\ell_4=0$.
But then $\iota_7=\ell_1$, and we can replace $d$ by
$d-c_\a^{(\b)}d'_0$ (cf.~(4.4) and (4.56)), so that $c_\a^{(\b)}$ becomes zero.
This proves that the assumption (4.39) does not holds. Thus we have (4.38).
\par
Now we prove
$$
d(t^{\v i})=0\mbox{ \ \ if \ \ }\v i=\v i_{_{\sc J_7}}.
\eqno(4.59)$$
 By Claim 1, we can suppose $|\v i|=n\ge2$. Assume
that we have proved (4.59) for $|\v i|<n$. Then $d([v,t^{\v
i}])=0$ for $v\in A$. {}From this, we obtain $d(t^{\v
i})\in\mbb{F}1_{\HH}$. Suppose $i_p>0$ for some $p\in I_7$. Then
$t^{\v i}=(i_{\ol p}+1)^{-1}[t^2_p,t^{\v i-\es_p+\es_{\ol p}}]$,
and $|\v i-\es_p+\es_{\ol p}|=n$, thus $d(t^{\v
i})\in[\mbb{F},t^{\v i-\es_p+\es_{\ol p}}]+[t^2_p,\mbb{F}]=0$.
Similarly, by replacing $d$ by $d-d'$ for some $d'\in\sum_{q\in
I_{2,3}\cup J_4}\mbb{F}\ptl_{t_q}$ (which does not affect the
results we have obtained so far), we can suppose
$$
d(x^{\a,\v i})=m_\a x^{\a,\v i}\;\;\;\for\;\;\;\a\in\G,\;\v i=\v
i_{_{\sc J_7}}, \eqno(4.60)$$ and
$$
d(t_p)=d(t^2_q)=0\;\;\;\for\;\;\;p\in I_{2,3}\cup J_4,\; q\in \ol
I_5\cup J_6. \eqno(4.61)$$ Note that $\HH$ is generated by
elements in (4.60) and (4.61),
 thus we obtain that (4.60) holds
for all $(\a,\v i)\in\G\times\JJ$. {}From this and (3.1), one can
easily deduce that
$$
\mu:\a\mapsto m_\a\mbox{ \ is a group homomorphism such that \ }
\mu\in\HOM\mbox{ \ if \ }\iota_7\ne\ell_1. \eqno(4.62)$$ Assume
that $\iota_7=\ell_1$. Then by (3.2) and (4.38), we have
$$
m_\a+m_\b=m_{\a+\b+\si_p}\mbox{ \ if \ }\a_p\b_{\ol p}\ne\a_{\ol p}\b_p
\mbox{ \ and }\a,\,\b\in\G,\;p\in I_{1,4}.
\eqno(4.63)$$
By (4.42), (4.45), (4.46), and by induction on $|i|+|j|$, one can prove
$$
m_{i\si_p+j\l_p}=0\;\;\;\;\for\;\;\;i,\,j\in\mb{Z},\;p\in I_{1,4}.
\eqno(4.64)$$
{}From this we want to prove
$$
m_\a=m_{\a+i\si_p+j\l_p}\;\;\;\for\;\;\;\a\in\G,\;i,j\in\mb{Z}.
\eqno(4.65)$$ By replacing $\a$ by some $\a+\si_p$ if necessary,
we can suppose $(\a_p,\a_{\ol p})\ne(0,0),(1,1)$. By (4.64) and by
$[x^\a,x^{-\si_p+j\l_p}]=(\a_p(i+j\l_p)-i\a_{\ol p})x^{\a+(i+1)\si_p+j\l_p}$,
we obtain $m_\a=m_{\a+i\si_p+j\l_p}$ if
$\a_p(i-1+je_p)\ne(i-1)\a_{\ol p}$, from this, one can deduce (4.65).
Now from (4.63) and (4.65), we obtain (4.62) again. Thus by
replacing $d$ by $d-d_\mu$, we have $d=0$. This proves Claim 4 and
also (4.10).
\par
To prove that (4.11) is a direct sum, suppose
$$
d=a'_0d'_0+\sum_{p\in\{0\}\cup J_1\cup\ol I_{2,3}\cup I_5}a_pd_p
+\sum_{q\in I_{2,3}\cup J_4}b_q\ptl_{t_q}
+d_\mu+\sum_{(0,0)\ne(\a,\v i)\in\G\times\JJ}c_{\a,\v
i}\ad_{x^{\a,\v i}}, \eqno(4.66)$$ is the $0$ derivation. Applying
$d$ to $A\cup\{1,t_p\,|\,p\in I_{2,3}\cup J_4\cup I_6\}$
(cf.~(4.36)), we obtain that all coefficients are zero except
$a'_0$. Thus $d=a'_0d'_0=0$. By (4.4), we obtain either $a'_0=0$ or
$d'_0=0$. Thus (4.11) is a direct sum. \qed\par \vs{5pt}
\par\ni
{\bf5. \ Second cohomology groups} \vs{-1pt}\par\ni In this
section, we shall determine the second cohomology groups of the
Hamiltonian Lie algebra $\HH=\HH(\v\ell,\G)$. It is well known
that all one-dimensional central extensions of a Lie algebra are
determined by the second cohomology group. Central extensions are
often used in the structure theory and the representation theory
of Kac-Moody algebras [K3]. Using central extension, we can
construct many infinite dimensional Lie algebras, such as affine
Lie algebras, infinite dimensional Heisenberg algebras, and
generalized Virasoro and super-Virasoro algebras, which have a
profound mathematical and physical background
(cf.~[K3,\,S1,\,SZ]). Since the cohomology groups are closely
related to the structures of Lie algebras, the computation of
cohomology groups seems to be important and interesting as well
(cf.~[J,\,LW,\,S1,\,S2,\,S3,\,SZ]).
\par
Recall that a {\it 2-cocycle} on $\HH$ is an $\mbb{F}$-bilinear function
$\psi:\HH\times \HH\rar
\mbb{F}$ satisfying the following conditions:
$$
\matrix{
\psi(v_1,v_2)=-\psi(v_2,v_1)\hfill&\mbox{(skew-symmetry)},
\vs{4pt}\hfill\cr
\psi([v_1,v_2],v_3)+\psi([v_2,v_3],v_1)+\psi([v_3,v_1],v_2)=0
\hfill&\mbox{(Jacobian identity)},
\hfill\cr}
\eqno\matrix{(5.1)\vs{4pt}\cr(5.2)\cr}$$
for $v_1,v_2,v_3\in \HH$.
Denote by $C^2(\HH,\mbb{F})$ the vector space of 2-cocycles on $\HH$.
For any $\mbb{F}$-linear function $f:\HH\rar\mbb{F}$, one can define a 2-cocycle
$\psi_f$ as follows
$$
\psi_f(v_1,v_2)=f([v_1,v_2])\;\;\;\for\;\;\;v_1,v_2\in \HH.
\eqno(5.3)$$
Such a 2-cocycle is called a {\it 2-coboundary} or a {\it trivial
2-cocycle} on $\HH$.
Denote by $B^2(\HH,\mbb{F})$ the
vector space of 2-coboundaries on $\HH$.
A 2-cocycle $\phi$ is said to be
{\it equivalent to} a 2-cocycle $\psi$ if $\phi-\psi$ is trivial. For
a 2-cocycle $\psi$, we denote by $[\psi]$ the equivalent class of $\psi$.
The quotient space
$$
H^2(\HH,\mbb{F})=C^2(\HH,\mbb{F})/B^2(\HH,\mbb{F})
=\{\mbox{the equivalent classes of 2-cocycles}\},
\eqno(5.4)$$
is called the {\it second cohomology group} of $\HH$.
\par\ni
{\bf Lemma 5.1}.
{\it If $\iota_7\ne\ell_1$, then $H^2(\HH,\mbb{F})=0$.}
\par\ni
{\it Proof.}
Let $\psi$ be a 2-cocycle.
Say $\ell_4\ne0$ (the proof is exactly similar if $\ell_i\ne0$ for $i\ne1,4$).
We fix $p\in I_4$.
Define a linear function $f$ by induction on $i_{\ol p}$
as follows:
$$
f(x^{\a,\v i})=\left\{\matrix{\a_{\ol p}^{-1}(\psi(t_p,x^{\a,\v
i}) -i_{\ol p}f(x^{\a,\v i-\es_{\ol p}}))\hfill&\mbox{if \
}\a_{\ol p}\ne0, \vs{4pt}\hfill\cr (i_{\ol
p}+1)^{-1}\psi(t_p,x^{\a,\v i+\es_{\ol p}})\hfill& \mbox{if \
}\a_{\ol p}=0, }\right. \eqno(5.5)$$ for $(\a,\v
i)\in\G\times\JJ$. Set $\phi=\psi-\psi_f$. Then (5.5) shows that
$$
\phi(t_p,x^{\a,\v i})=0\;\;\;\for\;\;\;(\a,\v i)\in\G\times\JJ.
\eqno(5.6)$$
Using Jacobian identity (5.2), we obtain
$$
0=\phi(t_p,[x^{\a,\v i},x^{\b,\v j}])=
(\a_{\ol p}+\b_{\ol p})\phi(x^{\a,\v i},x^{\b,\v j})
+i_{\ol p}\phi(x^{\a,\v i-\es_{\ol p}},x^{\b,\v j})
+j_{\ol p}\phi(x^{\a,\v i},x^{\b,\v j-\es_{\ol p}}),
\eqno(5.7)$$
for $(\a,\v i),\,(\b,\v j)\in\G\times\JJ$. If $\a_{\ol p}+\b_{\ol p}\ne0$,
by induction on $i_{\ol p}+j_{\ol p}$, we obtain
$\phi(x^{\a,\v i},x^{\b,\v j})=0$. On the other hand, if
$\a_{\ol p}+\b_{\ol p}=0$, then (5.7) gives
$$
\phi(x^{\a,\v i},x^{\b,\v j}) =-j_{\ol p}(i_{\ol
p}+1)^{-1}\phi(x^{\a,\v i+\es_{\ol p}}, x^{\b,\v j-\es_{\ol p}}),
\eqno(5.8)$$ and by induction on $j_{\ol p}$, we again have
$\phi(x^{\a,\v i},x^{\b,\v j})=0$. Thus $\phi=0$. \qed\par Assume
that $\iota_7=\ell_1$. Denote $\si=\sum_{p\in I_1}\si_p$, and we
use notation $\hom$ as in (4.9) (cf.~(4.6)). We construct
2-cocycles $\phi_p,\,\phi'_p,\,\phi_\mu$ for $p\in
I_1,\,\mu\in\hom$ as follows:
$$
\matrix{
\phi_p(x^\a,x^\b)=\a_p\d_{\a+\b,\si-\si_p},\vs{4pt}\hfill\cr
\phi'_p(x^\a,x^\b)=\a_{\ol
p}\d_{\a+\b,\si-\si_p},\vs{4pt}\hfill\cr
\phi_\mu(x^\a,x^\b)=\mu(\a)\d_{\a+\b,\si}, \hfill\cr}
\eqno\matrix{\hfill(5.9)\vs{4pt}\cr(5.10)\vs{4pt}\cr(5.11)\cr}$$
for $\a,\b\in\G$. It is straightforward to verify that they are
2-cocycles (cf.~[J]). {}From the proof of Theorem 5.2 below, one
can see why we construct such 2-cocycles.
\par\ni
{\bf Theorem 5.2}. {\it (1) $H^2(\HH,\mbb{F})=0$ if
$\iota_7\ne\ell_1$; (2) if $\iota_7=\ell_1$, then
$H^2(\HH,\mbb{F})$ is the vector space spanned by
$B=\{[\phi_p],[\phi'_p],[\phi_\mu]\,|\,p\in I_1,\mu\in\hom\}$.
Furthermore, for $a_p,b_p\in\mbb{F},\,\mu\in\hom$, we have
$$
\sum_{p\in I_1}(a_p[\phi_p]+b_p[\phi'_p])+[\phi_\mu]=0\ \Lra\
a_p=b_p=\mu=0. \eqno(5.12)$$ }
\par\ni
{\it Proof}. (1) follows from Lemma 5.1, while (2) follows from
[J]. We give a simple proof of (2) as follows.
\par
First we prove (5.12). Thus suppose
$$
\psi=\sum_{p\in I_1}(a_p\phi_p+b_p\phi'_p)+\phi_\mu+\psi_f,
\eqno(5.13)$$ is the zero 2-cocycle for some $a_p,b_p\in\mbb{F}$
and some linear function $f$. Then for $p\in I_1,\,\a\in\G$, by
applying $\psi$ to
$(x^{-\si_p},x^\si),\,(x^{\l_p},x^{\si-\l_p-\si_p}),\,
(x^\a,x^{\si-\a})$, we have
$$
\matrix{0=\psi(x^{-\si_p},x^\si)=-a_p-b_p, \vs{4pt}\hfill\cr
0=\psi(x^{\l_p},x^{\si-\l_p-\si_p})=e_pb_p, \vs{4pt}\hfill\cr\dis
0=\psi(x^\a,x^{\si-\a})=\mu(\a)+ \sum_{p\in I_1}(\a_p-\a_{\ol
p})f(x^{\si_p+\si}), \hfill\cr}
\eqno\matrix{(5.14)\vs{4pt}\hfill\cr(5.15)\vs{4pt}\cr(5.16)\cr}$$
(cf.~the definition of $\l_p$ in (3.97)). We obtain that
$a_p=b_p=0$ for $p\in I_1$ and by (4.9),
$$
\mu=\sum_{p\in I_1}c_p\mu_p\in
\hom\cap{\rm span}\{\mu_p\,|\,p\in I_1\}=\{0\},
\eqno(5.17)$$
where $c_p=-f(x^{\si_p+\si})\in\mbb{F}$. This proves (5.12).
\par
Now suppose $\psi$ is a 2-cocycle. We define a linear function $f$ as follows:
set $f(x^\si)=0$, and for $\a\in\G\bs\{\si\}$,
we define
$$
p_\a={\rm min}\{p\in I_1\,|\,(\a_p,\a_{\ol p})\ne(1,1)\},
\eqno(5.18)$$
and set
$$
f(x^\a)=\left\{\matrix{(\a_p-\a_{\ol
p})^{-1}\psi(x^{-\si_p},x^\a)\hfill& \mbox{if \ }\a_{\ol
p}\ne\a_p, \vs{4pt}\hfill\cr
e_p^{-1}(1-\a_p)^{-1}\psi(x^{\l_p},x^{\a-\si_p-\l_p})\hfill&
\mbox{if \ }\a_{\ol p}=\a_p\ne1, \hfill\cr}\right. \eqno(5.19)$$
for $p=p_\a$. Set
$$
\phi=\psi-\sum_{p\in I_1}(a_p\phi_p+b_p\phi'_p)-\psi_f,
\eqno(5.20)$$ where
$$
a_p=-\psi(x^{-\si_p},x^\si)-b_p,\,\;\;b_p=
e_p^{-1}\psi(x^{\l_p},x^{\si-\l_p-\si_p}), \eqno(5.21)$$
(cf.~(5.14) and (5.15)). Then one can prove
$$
\phi(x^{-\si_p},x^\a)=0\;\;\;\for\;\;\;p\in I_1,\;\a\in\G.
\eqno(5.22)$$ In fact, if $\a=\si$, it follows from (5.9), (5.10),
(5.20) and (5.21). Assume $\a\ne\si$. Let $p=p_\a$ and write
$$
x^\a=\left\{\matrix{(\a_p-\a_{\ol p})^{-1}[x^{-\si_p},x^\a]\hfill&
\mbox{if \ }\a_{\ol p}\ne\a_p, \vs{4pt}\hfill\cr
e_p^{-1}(1-\a_p)^{-1}[x^{\l_p},x^{\a-\si_p-\l_p}]\hfill&\mbox{if \
} \a_{\ol p}=\a_p\ne1, \hfill\cr}\right. \eqno(5.23)$$
(cf.~(5.19)), we can obtain (5.22) by the Jacobian identity (5.2).
{}From (5.22), by considering $\phi(x^{-\si_p},[x^\a,x^\b])$ and
by (5.2), we obtain
$$
\phi(x^\a,x^\b)=0\;\;\;\mbox{if}\;\;\;\a_p+\b_p\ne\a_{\ol p}+\b_{\ol p}\;\;\;
\mbox{for some}\;\;\;p\in I_1.
\eqno(5.24)$$
Now we want to prove
$$
\phi(x^{\l_q},x^\a)=0\;\;\;\for\;\;\;\a\in\G,\;q\in I_1.
\eqno(5.25)$$ By (5.24), we can suppose $\a_p=\a_{\ol p}$ if $p\ne
q$, and $\a_q=\a_{\ol q}+e_q$. If $\a=\si-\l_q-\si_q$, (5.25)
follows from (5.9), (5.10), (5.20) and (5.21). So suppose
$\a\ne\si-\l_q-\si_q$. Let $\a'=\a+\l_q+\si_q\ne\si$. If
$p=p_{\a'}\ne q$, then $\a'_p=\a'_{\ol p}=\a_p=\a_{\ol p}\ne1$,
and by writing
$x^\a=e_p^{-1}(1-\a_p)^{-1}[x^{\l_p},x^{\a-\si_p-\l_p}]$, we
obtain
$$
\phi(x^{\l_q},x^\a)=-e_p^{-1}(1-\a_p)^{-1}e_q\a_q\phi(x^{\l_p},x^{\a'-\si_p-\l_p})=0,
\eqno(5.26)$$ by (5.19). On the other hand, if $p=q$, we again
have $\phi(x^{\l_q},x^\a)=\phi(x^{\l_q},x^{\a'-\si_p-\l_p})=0$ by
(5.19).
\par
Now by (5.25) and
by writing $(\a_p-1)x^\a=-e_p^{-1}[x^{\l_p},x^{\a-\l_p-\si_p}]$, we obtain
$$
(\a_p-1)\phi(x^\a,x^\b)=-\b_p\phi(x^{\a-\l_p-\si_p},x^{\b+\l_p+\si_p}),
\eqno(5.27)$$
$$
(\a_p-2)\phi(x^{\a-\l_p-\si_p},x^{\b+\l_p+\si_p})=
-(\b_p+1)\phi(x^{\a-2\l_p-2\si_p},x^{\b+2\l_p+2\si_p}),
\eqno(5.28)$$ for $p\in I_1$,
 where (5.28) is obtained from
(5.27). Using (5.28) in (5.27) and by writing
$$
(3(\a_{\ol p}-\a_p)-2\a_pe_p)x^{\a-2\l_p-2\si_p}=[x^\a,x^{-2\l_p-3\si_p}],
\eqno(5.29)$$
we obtain
$$
\matrix{
(3(\a_{\ol p}-\a_p)-2\a_pe_p)(\a_p-1)(\a_p-2)\phi(x^\a,x^\b)
\vs{4pt}\hfill\cr\ \ \ \
=\b_p(\b_p+1)(\phi([x^\a,x^{\b+2\l_p+2\si_p}],x^{-2\l_p-3\si_p})
+\phi(x^\a,[x^{-2\l_p-3\si_p},x^{\b+2\l_p+2\si_p}])).
\hfill\cr}
\eqno(5.30)$$
We prove that
$$
\phi(x^\a,x^{-2\l_p-3\si_p})=0\;\;\;\for\;\;\;\a\in\G,\;p\in I_1.
\eqno(5.31)$$
By (5.24), we can suppose $\a_{\ol p}=\a_p+2e_p$. If $\a_p\ne1,2$, by
setting $\b=-2\l_p-3\si_p$ in (5.27) and (5,28), then the right-hand side
of (5.28) is zero by (5.22), and thus (5.31) holds. Suppose $\a_p=1,2$. Then
$\a_{\ol p}=1+2e_p$ or $2+2e_p$, thus we can
write $\a$ in the following form
$$
\a=\a'+\si_p+2\l_p\mbox{ \,or \,}\a'+2\si_p+2\l_p
\mbox{ for some }\a'\in\G
\mbox{ such that }(\a'_p,\a'_{\ol p})=(0,0).
\eqno(5.32)$$
We denote
$$
\matrix{
c_i=\phi(x^{\a'+i(\si_p+\l_p)-\si_p},x^{-i(\si_p+\l_p)-\si_p}),\hfill&
c'_i=\phi(x^{\a'+i(\si_p+\l_p)},x^{-i(\si_p+\l_p)-\si_p}),
\vs{4pt}\hfill\cr
d_i=\phi(x^{\a'+i(\si_p+\l_p)-2\si_p},x^{-i(\si_p+\l_p)}),\hfill&
d'_i=\phi(x^{\a'+i(\si_p+\l_p)-\si_p},x^{-i(\si_p+\l_p)}),
\hfill\cr}
\eqno\matrix{(5.33)\vs{4pt}\cr(5.34)\cr}$$
for $i\in\mb{Z}.$ By writing
$$
(i-2j)e_px^{-i(\si_p+\l_p)-\si_p}=
[x^{-j(\si_p+\l_p)-\si_p},x^{-(i-j)(\si_p+\l_p)-\si_p}]
\;\;\;\for\;\;\;j\in\mb{Z},
\eqno(5.35)$$
we obtain
$$
(i-2j)c_i=(i+j)c_{i-j}-(2i-j)c_j,\;\;\;
(i-2j)c'_i=i(c'_{i-j}-c'_j)\;\;\;\for\;\;\;i,j\in\mb{Z}.
\eqno(5.36)$$
By writing
$$
\matrix{ 2(j-i)e_px^{\a'+i(\si_p+\l_p)-\si_p}=
[x^{\a'+j(\si_p+\l_p)-2\si_p},x^{(i-j)(\si_p+\l_p)}],
\vs{4pt}\hfill\cr (j-i)e_px^{\a'+i(\si_p+\l_p)}=
[x^{\a'+j(\si_p+\l_p)-\si_p},x^{(i-j)(\si_p+\l_p)}], \hfill\cr}
\eqno\matrix{(5.37)\vs{4pt}\cr(5.38)\cr}$$ we obtain
$$
2(j-i)c_i=(2i+j)d_{j-i}-(i-j)d_j,\;\;
(j-i)c'_i=(i+j)d'_{j-i}+(i-j)d'_j\;\;\;\for\;\;\;i,j\in\mb{Z}.
\eqno(5.39)$$ Note that the system (5.36) has up to multiplicative
scalars unique solutions for $c_i,c'_i$, and we find that
$$
c_i=(i^3-i)c,\;\;\;c'_i=i^2c'\;\;\;\mbox{for \,$i\in\mb{Z}$\,\ and
some}\;\;\;c,c'\in\mbb{F},
\eqno(5.40)$$
are the only solutions. If we substitute
$j$ by $1$ and by $i+1$ in (5.39), we then obtain $c_i=c'_i=d_i=d'_i=0$
for all $i\in\mb{Z}$. This in particular proves (5.31) by (5.32)-(5.34).
\par
Now using (5.31) in (5.30), noting that
$\b_p-\b_{\ol p}=\a_{\ol p}-\a_p$ by (5.24), we deduce that
$$
(\b_p(\b_p+1)(3(\a_{\ol p}-\a_p)+2e_p(\b_p-1))
-(\a_p-1)(\a_p-2)(3(\a_{\ol p}-\a_p)-2\a_pe_p))\phi(x^\a,x^\b)=0.
\eqno(5.41)$$
As in the proof of (5.31), we can prove $\phi(x^\a,x^{2\l_p})=0$. Thus we can
replace $\l_p$ by $2\l_p$ in the above discussion, i.e., if we replace
$e_p$ by $2e_p$, (5.41) still holds. This forces
$$
\phi(x^\a,x^\b)=0\mbox{ \ or \ }\b_p-1=-\a_p\mbox{ \ for all \ }
p\in I_1,
\eqno(5.42)$$
i.e., if $\a+\b\ne\si$, then $\phi(x^\a,x^\b)=0$.
Thus we can suppose
$$
\phi(x^\a,x^\b)=m_\a\d_{\a+\b,\si}\mbox{ \ for \ }\a,\b\in\G
\mbox{ \,and some \,}m_\a\in\mbb{F}. \eqno(5.43)$$ As in the proof
of (5.31), we can prove $m_{i\si_p+j\l_p}=0$ for
$i,j\in\mb{Z},\,p\in I_1$. Then for any $\a,\b\in\G,\,p\in I_1$,
let $v_1=x^\a,v_2=x^\b, v_3=x^{\si-\a-\b-\si_p}$ in (5.2), one can
easily deduce that $\mu:\a\mapsto m_\a$ is a group homomorphism
$\mu:\G\to\mbb{F}$ such that $\mu(\si_p)=0$.
 Thus $\mu\in\HOM$ and
$\phi=\psi_\mu$. Furthermore, we can write $\mu=\nu+\l$ for
$\nu\in\hom,\,\l\in{\rm span} \{\mu_p\,|\,p\in I_1\}$ by (4.9).
Then $\phi=\psi_\nu+\psi_\l$. But from (5.16) and (5.17), one can
see that $\psi_\l$ corresponds to a trivial 2-cocycle, thus we can
suppose $\phi=\phi_\nu$. This proves Theorem 5.1.
\qed
\vs{5pt}\par\ni {\bf References} \vs{-1pt}\begin{description}
\item[{[DZ]}] D.~Dokovic and K.~Zhao,
Derivations, isomorphisms and second cohomology of generalized Block
algebras, Alg.~Colloq. {\bf 3} (1996), 245--272.

\item[{[F]}] R.~Farnsteiner, Derivations and central extensions of
  finitely generated graded Lie algebras, J.~Alg. {\bf118} (1988), 33--45.

\item[{[J]}] Y.~Jia, Second cohomology of generalized Cartan type $H$ Lie
algebras in characteristic 0, J.~Alg. {\bf204} (1998), 312--323.

\item[{[K1]}] V.~G.~Kac,
A description of filtered Lie algebras whose associated graded Lie algebras
are of Cartan types, Math.~of USSR-Izvestijia {\bf 8} (1974), 801--835.

\item[{[K2]}] V.~G.~Kac,
Classification of infinite-dimensional simple linearly compact Lie
superalgebras, Adv.~Math. {\bf 139} (1998), 1--55.

\item[{[K3]}] V.~G.~Kac, {\it Infinite Dimensional Lie Algebras}, 3rd ed.,
Cambridge Univ.~Press, 1990.

\item[{[LW]}] W.~Li, R.~L.~Wilson,
Central extensions of some Lie algebras,
Proc.~Amer. Math.~Soc. {\bf126} (1998), 2569--2577.

\item[{[O]}] J.~M.~Osborn,
New simple infinite-dimensional Lie algebras of characteristic 0,
J.~Alg. {\bf 185} (1996), 820--835.

\item[{[OZ]}] J.~M.~Osborn and K.~Zhao,
Generalized Poisson brackets and Lie algebras for type $H$ in
characteristic 0, Math.~Z. {\bf 230} (1999), 107--143.

\item[{[S1]}] Y.~Su, 2-Cocycles on the Lie algebras of generalized
differential operators, Comm.~Alg. {\bf 30} (2002), 763--782.

\item[{[S2]}] Y.~Su, On the low dimensional cohomology of Kac-Moody algebras
with coefficients in the complex field,
Shuxue Jinzhan (Chinese Adv.~in Math.) {\bf1989}, {\it18}, 346-351.

\item[{[S3]}] Y.~Su, 2-Cocycles on the Lie algebras of all
differential operators of several indeterminates.
(Chinese) Northeastern Math.~J. {\bf6} (1990), 365--368.

\item[{[SX]}] Y.~Su and X.~Xu, Central simple Poisson algebras, to appear.

\item[{[SXZ]}] Y.~Su, X.~Xu and H.~Zhang,
Derivation-simple algebras and the structures of Lie algebras of Witt type,
J.~Alg. {\bf 233} (2000), 642--662.

\item[{[SZ]}] Y.~Su, K.~Zhao, Second cohomology group of
generalized Witt type Lie algebras and certain representations,
Comm.~Alg. {\bf 30} (2002), 3285--3309.

\item[{[X1]}] X.~Xu, Generalizations of Block algebras,
  Manuscripta Math. {\bf 100} (1999), 489--518.

\item[{[X2]}] X.~Xu,
New generalized simple Lie algebras of Cartan type over a field
with characteristic 0, J.~Alg. {\bf 224} (2000), 23--58.

\item[{[Z]}] H.~Zhang, The representations of the coordinate ring
 of the quantum symplectic space, J.~Pure Appl.~Algebra {\bf150} (2000), 95--106.

\item[{[Zh]}] K.~Zhao,
A class of infinite dimensional simple Lie algebras,
J.~London Math.~Soc.~(2), {\bf62} (2000), 71--84.
\end{description}
\end{document}